\documentclass[leqno,11pt]{amsart}

\usepackage{amssymb}
\usepackage{amsmath}
\usepackage[usenames, dvipsnames]{color}
\usepackage{enumitem}
\usepackage{amsthm}
\usepackage{mathrsfs}
\usepackage{amsrefs}
\numberwithin{equation}{section}

\usepackage{marginnote}
\usepackage{todonotes}
\usepackage{hyperref}
\usepackage{esint}
\usepackage{verbatim}
\usepackage{tikz}
\usepackage{graphicx}
\usetikzlibrary{positioning, arrows.meta}

\usepackage{geometry} 

\usepackage{etoolbox}

\theoremstyle{plain}
\newtheorem{theorem}{Theorem}[section]
\newtheorem{lemma}[theorem]{Lemma}
\newtheorem{corollary}[theorem]{Corollary}
\newtheorem{proposition}[theorem]{Proposition}

\newtheorem{definition}[theorem]{Definition}

\newtheorem{remark}[theorem]{Remark}

\newcommand{\R}{\mathbb{R}}
\newcommand{\N}{\mathbb{N}}

\newcommand{\bA}{\mathbf{A}}

\newcommand{\M}{\mathcal{M}}

\newcommand{\W}{\mathcal{W}}
\newcommand{\hW}{\hat{\mathcal{W}}}
\newcommand{\mV}{\mathrm{V}}

\makeatletter
\@namedef{subjclassname@2020}{\textup{2020} Mathematics Subject Classification}
\makeatother

\newcommand{\wei}[1]{\langle #1 \rangle}

\def\XXint#1#2#3{{\setbox0=\hbox{$#1{#2#3}{\int}$ }
\vcenter{\hbox{$#2#3$ }}\kern-.6\wd0}}

\title[Equations with singular-degenerate coefficients, Conormal problem]{non-homogeneous parabolic conormal boundary problem with singular-degenerate coefficients}

\author[J. Fang]{Junyuan Fang}
\address[J. Fang]{Department of Mathematics, Louisiana State University, 125B Prescott Hall,
Baton Rouge, LA 70803-4918}
\email{jfang8@lsu.edu}

\author[T. Phan]{Tuoc Phan}
\address[T. Phan]{Department of Mathematics, University of Tennessee, 227 Ayres Hall,
1403 Circle Drive, Knoxville, TN 37996-1320}
\email{tphan2@utk.edu}


\begin{document}

\begin{abstract}  We study non-homogeneous conormal boundary value problems for a class of divergence-form parabolic equations whose volumetric heat-capacity coefficients belong to a Muckenhoupt class of weights. Consequently, the equations exhibit singular, degenerate, or simultaneously singular-degenerate behavior. Under smallness conditions on the partial mean oscillation of the diffusion coefficients and the weighted mean oscillation of the volumetric heat-capacity coefficient, we establish regularity estimates, existence, and uniqueness results in appropriate weighted Sobolev spaces. To obtain these results, we employ the level-set argument introduced by Caffarelli and Peral. We also establish weighted parabolic cylinders, weighted inequalities, and a weighted Aubin-Lions compactness theorem, which are essential tools in our analysis.   
\end{abstract}
\maketitle

\section{Introduction} 
\subsection{Problem settings and motivations} \label{mov-sec} Let $\Omega_T = \Omega \times (0, T]$, where  $T \in (0,\infty)$ is a fixed number, and $\Omega \subset \mathbb{R}^n$ is a non-empty open connected bounded set with $C^1$-boundary $\partial \Omega$.  We study the following non-homogeneous conormal boundary value problem for the class of parabolic equations with singular-degenerate volumetric heat-capacity coefficients
\begin{equation} \label{main-eqn}
\left\{
\begin{array}{ccccl}
\beta(x) u_t - \textup{div}(\bA(x,t) \nabla u) &= & \textup{div}(F)  &\quad  &\text{in} \quad \Omega_T, \\
\wei{\bA(x,t) \nabla u +F, \vec{\mathbf{n}}} & =  &g &\quad &\text{on} \quad \partial\Omega \times (0,T), \\
u(\cdot,0) &= & u_0(\cdot) &\quad &\text{in} \quad \Omega.
\end{array}\right.
\end{equation}
Here in \eqref{main-eqn}, $F: \Omega_T \rightarrow \mathbb{R}^n, u_0: \Omega \rightarrow \mathbb{R}$,  and $g: \partial \Omega \times (0, T) \rightarrow \mathbb{R}$ are given data which are assumed to be in suitable spaces, and $\vec{\mathbf{n}}$ is the unit outward normal vector on $\partial \Omega$. The diffusion coefficient $\bA: \Omega_T \rightarrow \mathbb{R}^{n\times n}$ is a measurable matrix-valued function satisfying the following ellipticity and boundedness conditions: there is a constant $\nu \in (0,1)$ such that
\begin{equation} \label{ellip-cond}
\nu |\xi|^2 \leq \wei{\bA(x,t) \xi, \xi} \quad \text{and} \quad |\bA(x,t)| \leq \nu^{-1}, 
\end{equation}
for all $(x,t) \in \Omega_T$ and for all $\xi=(\xi_1, \xi_2,\ldots, \xi_n) \in \mathbb{R}^n$, with $|\bA|=\sup_{|\xi|=1}|\bA\xi|$ denoting the Euclidean operator norm. Moreover, the volumetric heat-capacity coefficient $\beta: \mathbb{R}^n\to[0,\infty)$ is measurable such that $\frac{1}{\beta}$ is in a suitable Muckenhoupt class of weights specified in \eqref{beta-cond} below. As such, the coefficients in \eqref{main-eqn} can be singular, degenerate or simultaneously singular and degenerate. Weak solutions to the problem \eqref{main-eqn} are understood in the sense defined in \eqref{weak-sol-def} below. When $\beta\equiv1$, the equation \eqref{main-eqn} reduces to the classical second-order parabolic equation in divergence form with inhomogeneous conormal boundary condition.

\smallskip
The singular-degenerate model \eqref{main-eqn}, characterized by the heat-capacity coefficient $\beta$, arises naturally in the study of heterogeneous diffusion and phase-transition phenomena, where the effective heat-storage capacity or thermal inertia varies significantly across the spatial domain. For instance, in media containing vacuum regions or low-density gases, both the heat capacity and thermal conductivity may vanish locally; see \cite{DIT, Li-Xin, PSW}. Similarly, in heat conduction through composite materials and in phase-transition processes, the heat capacity may undergo abrupt changes across different materials or become singular near phase interfaces, leading to singular or degenerate parabolic equations. We refer the reader to \cite{AS, FKRH, PVPF} and the references therein for further discussions. The non-homogeneous conormal boundary condition considered in this paper naturally models prescribed boundary fluxes and transmission phenomena across material interfaces.

\smallskip
The main objective of this paper is to establish the existence, well-posedness, and regularity theory for problem \eqref{main-eqn} in appropriate weighted Sobolev spaces under the assumptions that the initial datum $u_0$ belongs to $W^{1,p}(\Omega)$, the source term $F$ is in $L^p(\Omega_T)^n$, and the boundary datum $g$ belongs to $L^p(\partial\Omega\times(0,T))$, where $p\in(1,\infty)$, and $\frac{1}{\beta}$ belongs to a suitable Muckenhoupt class of weights specified in \eqref{beta-cond}. A significant difficulty is that no time regularity is assumed for $g$; in particular, $g_t$ need not exist, and hence the standard reduction $u\mapsto u-g$ is not available.
Another challenge arises from assumption \eqref{beta-cond}, which allows the heat-capacity coefficient $\beta:\mathbb{R}^n\to[0,\infty)$ to vanish or become unbounded at points in $\overline{\Omega}$. Consequently, equation \eqref{main-eqn} may exhibit both singular and degenerate behavior which requires subtle analysis and estimates.

\smallskip
From an analytical perspective, incorporating a Muckenhoupt weight as in \eqref{beta-cond} into the coefficient structure of \eqref{main-eqn} provides a mathematically rigorous framework capable of simultaneously capturing spatial singularities and degeneracies. This weighted approach has played a crucial role in establishing a priori estimates, weighted Sobolev embeddings, and regularity theories for linear and nonlinear elliptic and parabolic equations. This approach was initiated in \cite{Fabes, Fabes-1} where H\"{o}lder regularity estimates are proved, and also in the recent work such as \cite{CFP, DPS, DP, DPT1, DPT2, FP1, FP2}, which contain various results on regularity in weighted and mixed-norm weighted Sobolev spaces and on Krylov-Safonov Harnack inequalities. This paper extends the line of research initiated in \cite{FP1} by transitioning from the homogeneous Dirichlet problem to the non-homogeneous conormal boundary value problem for \eqref{main-eqn}.
\subsection{Main results} Let us begin by introducing several notations and definitions.  Throughout the paper, we assume that the weight $\beta$ in \eqref{main-eqn} satisfies the condition
\begin{equation} \label{beta-cond}
\frac{1}{\beta}\in A_{1+\frac{2}{n_0}} \quad \text{with} \quad \big[\frac{1}{\beta}\big]_{A_{1+\frac{2}{n_0}}} \leq M_0, \quad \text{where}\quad n_0 = \max\{n, 2\}, 
\end{equation}
for some given $M_0 \geq 1$. See Definition \ref{A-p-def} below for the definition of the $A_{p}$ class of Muckenhoupt weights with $p \in (1, \infty)$. We note that when $n=1$ or $n=2$, condition \eqref{beta-cond} is equivalent to
\[
\beta \in A_2 \quad \text{with} \quad [\beta]_{A_2} \leq M_0.
\]
In addition, for $r>0$ and $x_0 \in \mathbb{R}^n$, $B_r(x_0)$ denotes the ball in $\mathbb{R}^n$ of radius $r$ centered at $x_0$. As introduced in \cite{FP1}, for $z_0 = (x_0, t_0) \in \mathbb{R}^n \times \mathbb{R}$, 
we define the weighted parabolic cylinder $Q_{r,\beta}(z_0)$ with radius $r>0$ by
\begin{equation} \label{cylinder-def}
Q_{r, \beta}(z_0) =B_r(x_0)  \times \Gamma_{\beta, z_0}(r), \quad \text{where}   \quad \Gamma_{\beta, z_0} (r) =(t_0 -r^2\mathcal{E}_{\beta}(x_0, r), t_0],
\end{equation}
and $\mathcal{E}_{\beta}(x_0, r)$ is defined by
\begin{equation} \label{t-height}
\mathcal{E}_{\beta}(x_0, r)= \Big( \fint_{B_{r}(x_0)} \beta(x)^{\frac{n_0}{2}} dx \Big)^{\frac{2}{n_0}},\quad \text{with}\quad n_0=\max\{n, 2\}.
\end{equation}
Due to \eqref{beta-cond}, we see that $\mathcal{E}_{\beta}(x_0, r)<\infty$. 

\smallskip
For a locally integrable function $f$ defined on $\Omega\subset \R^{n}$ and for a measurable subset $E\subset \Omega$ with a positive Lebesgue measure, we denote 
\[
(f)_{E}=\fint_{E}f(x)dx = \frac{1}{|E|} \int_{E} f(x) dx,
\]
where $|E|$ denotes the Lebesgue measure of $E$. As introduced in \cite{KK2} and in \cite{CMP}, the partial mean oscillation of $\bA$ and the weighted mean oscillation of $\beta$ are respectively defined by 
\begin{equation} \label{BMO-defi}
\begin{split}
& [[\bA]]_{\textup{BMO}}(\Omega_T, R_0)   = \sup_{\substack{z_0 =(x_0, t_0) \in \Omega_T \\ \rho \in (0, R_0)}}\Big(\fint_{Q_{\rho, \beta}(z_0) \cap \Omega_T} |\bA(x,t) - (\bA)_{B_\rho(x_0) \cap \Omega} (t)|^2 dxdt\Big)^{1/2}, \\[4pt]
& [[\beta]]_{\textup{BMO}}(\Omega, R_0) = \sup_{\substack{\rho \in (0, R_0)\\ x_0 \in \Omega}}\Big(\frac{1}{\beta(B_{\rho}(x_0))}\int_{B_\rho(x_0)} | \beta(x) - (\beta)_{B_\rho(x_0)}|^2 \beta(x)^{-1}dx\Big)^{1/2}.
\end{split}
\end{equation}
Here, in \eqref{BMO-defi}, $\beta(B_{\rho}(x_0))$ is the measure of the ball $B_{\rho}(x_0)$ with respect to the measure $\beta(x)dx$:
\[
\beta(B_{\rho}(x_0))=\int_{B_{\rho}(x_0)}\beta(x)dx.
\] 

\smallskip
Now, we refer the reader to Definition \ref{W-space} below for the definition of the weighted Sobolev space $\W^{1,p}(\Omega_T, \beta)$. A function $u \in \W^{1,p}(\Omega_T, \beta)$ is said to be a weak solution to \eqref{main-eqn}  if 
\begin{align} \notag
& -\int_{\Omega_T} \beta(x) u(x,t) \varphi_t(x,t) dx dt + \int_{\Omega_T} \wei{\bA(x,t) \nabla u(x,t), \nabla \varphi(x,t)} dx dt \\ \label{weak-sol-def}
& + \int_{\Omega_T} \wei{F(x,t), \nabla \varphi(x,t)} dx dt =  \int_{\Omega} \beta(x)u_0(x) \varphi(x, 0) dx + \int_{\partial \Omega \times (0,T)} g(x,t) \varphi(x,t) \,dSdt,
\end{align}
for every $\varphi \in C^\infty(\overline{\Omega} \times [0, T])$ that vanishes near $T$. The following theorem is the first main result of the paper.
\begin{theorem} \label{main-theorem} Let $p \in (1, \infty)$, $\nu \in (0,1)$, and $M_0\geq 1$. Then, there exists a sufficiently small constant $\delta = \delta (n, \nu, p, M_0)>0$ such that the following assertions hold. Suppose that  \eqref{ellip-cond} and \eqref{beta-cond} hold, $\partial \Omega \in C^1$, and
\begin{equation} \label{small-ness}
[[\bA]]_{\textup{BMO}}(\Omega_T, R_0)  + [[\beta]]_{\textup{BMO}}(\Omega, R_0) < \delta \qquad \text{for some} \quad R_0 \in (0,1).
\end{equation}
Then, for every $F \in L^p(\Omega_T)^n$ and $g\in L^p(\partial \Omega\times (0,T))$ with $T>0$, and for $u_0 \in W^{1,p}(\Omega)$, there exists a unique weak solution $u \in \W^{1,p}(\Omega_T, \beta)$ solving \eqref{main-eqn}. Moreover, 
\begin{equation} \label{main-thm-est}
\|u\|_{\W^{1,p}(\Omega_T,\beta)} \leq N \Big(\lambda_0^{1/p}\|u_0\|_{W^{1,p}(\Omega)} +\|F\|_{L^p(\Omega_T)} +\|g\|_{L^p(\partial \Omega\times (0,T))}  \Big),
\end{equation}
where $N = N(n, \nu, p, M_0, R_0, \Omega, T/\lambda_0) >0$ with
\begin{equation} \label{lambda-zero-def}
\lambda_0 = \left(\int_{\Omega} \beta(x)^{\frac{n_0}{2}} dx\right)^{\frac{2}{n_0}} \quad \text{for} \quad n_0 = \max\{2, n\}.
\end{equation}
\end{theorem}

\smallskip
Besides Theorem \ref{main-theorem}, we also investigate local boundary regularity estimates of solutions to the class of equations \eqref{main-eqn}. To this end, for $r>0$ and $z_0=(x_0,t_0)\in\overline{\R^n_+}\times\R$, we denote
\[
T_r(x_0) = B_r(x_0) \cap \{x_n = 0\} \quad \text{and}\quad Q_{r,\beta}^+(z_0)= B_r^+(x_0) \times \Gamma_{\beta, z_0}(r),
\]
where $B_r^+(x_0) = B_r(x_0) \cap \mathbb{R}^n_+$ with $\mathbb{R}^n_+ = \mathbb{R}^{n-1} \times (0, \infty)$, and $\Gamma_{\beta, z_0}(r)$ is defined in \eqref{cylinder-def}. When $x_0=0$ and $z_0=(0,0)$, we omit them from the notation and simply write
\[
T_r = T_r(0), \ B_r = B_r(0), \ B_r^+ = B_r^+(0), \ Q_{r, \beta}^+ = Q_{r, \beta}^+(0).
\]
Without loss of generality, we consider the following equations with non-homogeneous conormal boundary conditions on the flat boundary part
\begin{equation} \label{bdry eqn}
\left\{ 
\begin{aligned}
\beta(x) u_t - \textup{div}(\bA \nabla u) &= \textup{div}(F) \quad &&\text{in} \quad Q_{2\Lambda, \beta}^+, \\[4pt]
(\bA \nabla u + F)\cdot \vec{\mathbf{n}} &= g \quad &&\text{on}\quad  T_{2\Lambda}\times \Gamma_{\beta}(2\Lambda),
\end{aligned} \right.
\end{equation}
where $\vec{\mathbf{n}} = (0, \ldots, 0, -1) \in \mathbb{R}^n$, $\Lambda=\Lambda(n, M_0)>1$ is a number defined in Lemma \ref{quasi-metric lemma} below. The formal definition of weak solutions $u \in \hW^{1,p}(Q_{2\Lambda,\beta}^+,\beta)$ to \eqref{bdry eqn} is given in Definition~\ref{boundary-weak-def} below, in which the spaces $\hW^{1,p}(Q_{2\Lambda,\beta}^+,\beta)$ and $\mV_p(B_2^+)$ are defined in Definition \ref{space-hW1p}. The following theorem on regularity estimates of weak solutions to \eqref{bdry eqn} is another main result of the paper.  
\begin{theorem} \label{inter-theorem-bdry}
Let $\nu\in(0,1)$, $M_0\geq1$, and $p\geq 2$ be fixed. There exist constants $\delta_1=\delta_1(n,\nu,p,M_0)>0$ sufficiently small and $N=N(n,\nu,p,M_0)>0$ such that the following holds. Suppose that \eqref{beta-cond} holds, $\bA$ satisfies \eqref{ellip-cond} in $Q_{2\Lambda,\beta}^+$, and
\[
[[\bA]]_{\textup{BMO}}(Q_{1, \beta}^+, 1) + [[\beta]]_{\textup{BMO}}(B_{1}^+, 1) \leq \delta_1.
\]
Then for every weak solution $u\in\hW^{1,2}(Q_{2\Lambda,\beta}^+,\beta)$ of \eqref{bdry eqn} with $F\in L^p(Q_{2\Lambda,\beta}^+)^n$ and $g\in L^p(T_{2\Lambda}\times\Gamma_\beta(2\Lambda))$, it satisfies
\begin{align*}
&\|\nabla u\|_{L^p(Q_{1,\beta}^+)}
+\|\beta u_t\|_{L^p(\Gamma_\beta(1), \mV_{p'}(B_1^+)^{*})}\\
&\leq N\left(|Q_{1,\beta}^+|^{1/p-1/2}\|u\|_{L^2(Q_{2\Lambda,\beta}^+)}
+\|F\|_{L^p(Q_{2\Lambda,\beta}^+)}
+\|g\|_{L^p(T_{2\Lambda}\times\Gamma_\beta(2\Lambda))}\right),
\end{align*}
where $p'=p/(p-1)$, and $\mV_{p'}(B_1^+)^{*}$ is the dual of the space $\mV_{p'}(B_1^+)$, defined as in \eqref{space-Vp}.
\end{theorem}
\noindent We remark that the choice of $Q_{2\Lambda, \beta}^+$ in the theorem is for algebraic convenience in the proof. The general case for cylinders with different sizes can be derived using the dilation in Lemma~\ref{scaling-lemma} below.

\smallskip
Similar to Theorem \ref{inter-theorem-bdry}, we also investigate the regularity estimates for weak solutions to \eqref{main-eqn} locally in time and globally in the spatial domain. Specifically, we consider
\begin{equation}\label{main-eqn-2}
\left\{
\begin{aligned}
\beta(x) u_t - \textup{div}(\bA \nabla u) &= \textup{div}(F)\quad &&\text{in} \quad \Omega\times (0,T),\\
(\bA\nabla u+F)\cdot \vec{\mathbf{n}} &= g \quad &&\text{on} \quad \partial \Omega \times (0,T).
\end{aligned}\right.
\end{equation}
Note that in comparison with \eqref{main-eqn}, problem \eqref{main-eqn-2} imposes no initial condition. A formal definition of weak solutions to \eqref{main-eqn-2} is stated in Definition~\ref{interior-weak-def} below. The following result is our last main result of the paper.

\begin{theorem}\label{main-thm-2}Let $\nu \in (0,1)$, $M_0 \geq 1$, and $p \in [2, \infty)$. There exists a sufficiently small constant $\bar{\delta} = \bar{\delta}(n, \nu, p, M_0) \in (0,1)$ such that the following assertions hold. Suppose that  \eqref{ellip-cond} and \eqref{beta-cond} hold, $\partial \Omega \in C^1$, and
\begin{equation} \label{small-ness-inter-theorem}
[[\bA]]_{\textup{BMO}}(\Omega_T, R_0)  + [[\beta]]_{\textup{BMO}}(\Omega, R_0) < \bar{\delta} \qquad \text{for some} \quad R_0 \in (0,1).
\end{equation}
Then, for every weak solution $u \in \W^{1,2} (\Omega_T, \beta)$ to \eqref{main-eqn-2} with $F \in L^p(\Omega_T)^n$ and $g\in L^p(\partial \Omega\times (0, T))$, and for every $T_0 \in (0, T)$, it holds that $u \in \W^{1,p}(\Omega\times (T_0, T),\beta)$. Moreover, there exists a constant $N=N(n,\nu,p,M_0,R_0,\Omega,T_0/\lambda_0,T/\lambda_0)>0$ such that
\begin{equation} \label{global-time-cut-est}
\|u\|_{\W^{1,p}(\Omega\times (T_0, T),\beta)}\leq N\left(\lambda_0^{1/p-1/2}\|u\|_{L^2(\Omega_T)} + \|F\|_{L^p(\Omega_T)}+\|g\|_{L^{p}(\partial \Omega\times (0, T))} \right),
\end{equation}
where $\lambda_0>0$ is defined in \eqref{lambda-zero-def}.
\end{theorem}   
We remark that it is possible to relax the boundary smoothness condition $\partial\Omega \in C^1$ to some flatness conditions such as Lipschitz domains with sufficiently small Lipschitz constants and  Reifenberg-flat domains as those considered in \cites{CMP, B2, BW1, HNP}. In addition, extending the results of the paper to the weighted mixed norm setting following the method developed in \cite{DK, DP-0, DPT2} is possible as well. However, we do not pursue these extensions to avoid additional technical difficulties.

\subsection{\texorpdfstring{
Literature review and state of the art\unskip\unskip}{Literature review and state of the art}}
  As briefly discussed in Section \ref{mov-sec},  the study of equations with singular-degenerate coefficients governed by Muckenhoupt weights was initiated by Fabes, Kenig, and Serapioni in their seminal works \cite{Fabes, Fabes-1}. They essentially considered the following linear elliptic equations in divergence form
\begin{equation} \label{bar-A-wei}
\textup{div}(\mathbf{\bar{A}}(x)\nabla u)=0 \quad \text{in} \quad B_1,
\end{equation}
where the coefficient matrix $\mathbf{\bar{A}}(x)$ satisfies the weighted ellipticity condition
\begin{equation} \label{FKS-cond}
\nu \mu(x)|\xi|^2
\leq
\wei{\bar{\bA}(x)\xi,\xi}
\quad\text{and}\quad
|\bar{\bA}(x)|
\leq
\nu^{-1}\mu(x),
\qquad
\forall\,\xi\in\mathbb{R}^n,
\end{equation}
with a given $\mu\in A_2$. By establishing weighted Sobolev inequalities and adapting Moser's iteration technique, they proved an interior Harnack inequality and consequently obtained the local H\"{o}lder continuity of weak solutions. This foundational framework has since become a cornerstone of the regularity theory for singular-degenerate equations. It is important to note that there is a parallel line of research, initiated by Trudinger in his classical works \cite{Tru-1, Tru-2} and further developed in the recent breakthrough paper \cite{BS}. In these works, the smallest and largest eigenvalues of $\bar{\mathbf{A}}$ satisfy suitable integrability conditions instead of the structural assumptions imposed by Muckenhoupt weights, and H\"{o}lder regularity estimates were established. In addition, the parabolic counterpart of \eqref{bar-A-wei} in the framework of \cite{Fabes, Fabes-1} naturally includes space-time operators in which the temporal and spatial components are required to be balanced through singular or degenerate Muckenhoupt weights. We refer, for example, to \cite{Chia-1, Chia-2, Chia-3}, where Harnack inequalities and H\"{o}lder regularity estimates were established under additional structural assumptions on the underlying weights.

\smallskip
Moving beyond qualitative H\"{o}lder continuity, the development of a Calder\'{o}n-Zygmund theory for linear equations in divergence form with singular-degenerate coefficients requires a sophisticated adaptation of classical harmonic analysis to weighted function spaces. Consider the degenerate elliptic equation or system
\begin{equation} \label{CMP-eqn}
\textup{div}(\bar{\mathbf{A}}(x)\nabla u)=\textup{div}(F) \quad \text{in} \quad \Omega,
\end{equation}
where $\bar{\mathbf{A}}(x)$ satisfies the weighted ellipticity condition associated with an $A_2$ weight as in \eqref{FKS-cond}. As in the uniformly elliptic theory, higher integrability of the gradient cannot be expected without additional regularity assumptions on the coefficients. In particular, the coefficient matrix must exhibit sufficiently small oscillation in a weighted sense in order to obtain higher weighted Lebesgue estimates for $\nabla u$. Under suitable smallness assumptions on the weighted mean oscillation of $\bar{\mathbf{A}}$, a recent work \cite{CMP} established weighted Calder\'{o}n-Zygmund estimates in $W^{1,p}$ for elliptic equations. We also refer to \cite{BDGP}, where analogous results are obtained in the setting of matrix weights using a different weighted framework, and to \cite{DPS}, which develops Calder\'{o}n-Zygmund estimates for quasilinear elliptic equations with weights given by powers of the distance to the boundary that may not be in $A_2$.

\smallskip
On the other hand, for parabolic equations, the maximal $L^p$-regularity framework naturally requires simultaneous control of the time derivative $u_t$ and the gradient $\nabla u$ in weighted spaces consistent with parabolic scaling. Such estimates yield existence, uniqueness, and optimal regularity under minimal assumptions on the coefficients and the Muckenhoupt weights. Nevertheless, the available theory remains largely confined to highly structured weights, most notably powers of the distance to the boundary. A representative model is
\begin{equation}\label{DP-eqn}
\mu(x_n)u_t-\textup{div}\bigl(\mu(x_n)\mathbf{A}(x,t)\nabla u\bigr)
=\textup{div}\bigl(\mu(x_n)F\bigr),
\end{equation}
where $x=(x',x_n)\in\mathbb{R}^{n-1}\times(0,\infty)$, $\mu(x_n)\sim x_n^\alpha$ for some $\alpha\in\mathbb{R}$, and $\mathbf{A}$ is uniformly elliptic and bounded as in \eqref{ellip-cond}. Here, the singularity or degeneracy is localized solely along the boundary ${x_n=0}$. Well-posedness and mixed-norm weighted Sobolev estimates for equations of the form \eqref{DP-eqn}, as well as their non-divergence counterparts, have recently been established in the series of papers \cite{DP, DP-0, DP1, MNS-1} under various assumptions on $\alpha$ dictated by the boundary conditions.

\smallskip

In contrast to the balanced structure of \eqref{DP-eqn}, placing the weight exclusively in either the time derivative or the diffusion operator leads to fundamentally different analytical challenges. In such mismatched settings, even identifying an appropriate functional framework becomes highly nontrivial, and establishing a Sobolev regularity theory requires ideas beyond the classical weighted theory. Indeed, the pioneering work \cite{Chia-1} demonstrates that the standard H\"{o}lder regularity theory may fail for this class of equations. This fundamental difficulty has motivated a series of recent investigations, including \cite{DPT2, FP1, Phan-1} for divergence-form equations and \cite{DPT1, CFP, FP2} for their non-divergence counterparts. In particular, \cite{FP1} establishes well-posedness and weighted Sobolev regularity estimates for the homogeneous Dirichlet boundary value problem associated with equations having singular-degenerate volumetric heat-capacity coefficients of the form given in \eqref{main-eqn}:
\begin{equation}\label{FP-2.eqn}
\beta(x)u_t-\textup{div}\bigl(\mathbf{A}(x,t)\nabla u\bigr)
=\textup{div}(F),
\end{equation}
where the coefficient matrix $\mathbf{A}$ satisfies \eqref{ellip-cond}, while the weight $\beta$ satisfies \eqref{beta-cond}. The results established in this paper, Theorems \ref{main-theorem}--\ref{main-thm-2}, therefore provide the conormal counterpart of those in \cite{FP1}. 
See also the recent development in \cite{Phan-1} for similar results for homogeneous Dirichlet boundary value problems for the corresponding parabolic equations of \eqref{CMP-eqn} in which the singularity and degeneracy arise in the diffusion coefficients as \eqref{FKS-cond}, but not in the volumetric heat-capacity coefficients as those in \eqref{main-eqn} or \eqref{FP-2.eqn}. In addition, for the special case $\beta(x)=x_n^\alpha$ with
$x=(x',x_n)\in\mathbb{R}^{n-1}\times(0,\infty)$ and $\alpha\in(-2,0)$, analogous mixed-norm weighted Sobolev estimates for \eqref{FP-2.eqn} were obtained in \cite{DPT2} for divergence form equations and in \cite{DPT1} for the corresponding equations in non-divergence form, both with homogeneous Dirichlet boundary conditions. We also refer to \cite{CFP, FP2}, where Krylov-Safonov Harnack inequalities, H\"{o}lder regularity estimates, and Lin's $W^{2,\varepsilon}$ estimates have been established for the corresponding  non-divergence form of \eqref{FP-2.eqn} under different conditions on $\beta$.
\subsection{Main difficulties and our approaches}  
We apply the level set method introduced by Caffarelli and Peral in \cite{CP}, using a covering lemma and a freezing coefficient argument. Note that an alternative approach based on the Fefferman sharp function introduced by N. V. Krylov in \cite{Krylov-book} is possible; see \cite{KK2, DK, DP-0}. To implement the level set method, at each small parabolic cylinder with spatial base $B_r(x_0)$, we locally compare the solution $u$ of \eqref{main-eqn} with a solution $v$ of a reference equation with frozen coefficients:
\begin{equation}\label{reference-eqn}
(\beta)_{B_r(x_0)}v_t-\textup{div}\left((\bA)_{B_r(x_0)\cap\Omega}(t)\nabla v\right)=0.
\end{equation}
However, unlike the unweighted case $\beta\equiv1$, the equation for $u-v$ contains the additional term $(\beta(x)-(\beta)_{B_r(x_0)})v_t$, which is not directly controlled by the natural energy bounds. This prevents a straightforward energy estimate for $\nabla(u-v)$.

\smallskip
To overcome this obstruction, we use a compactness and contradiction argument in energy spaces that vary with $\beta$. The weighted Aubin--Lions theorem (Theorem~\ref{compactness-lemma} below), combined with weighted interpolation and localized energy identities, yields the gradient approximation established in Proposition~\ref{L2-comparision-bdry}. As in \cite{Phan-2}, we handle the boundary datum $g$ directly through trace inequalities, without requiring time regularity of $g$. Note also that it is possible to relax the integrability condition on $g$, but we did not pursue this extension to avoid further technical difficulties. Finally, the boundary Lipschitz estimate (Proposition~\ref{Bdry-Lipschitz-constant}) and a covering argument adapted to the weighted cylinders allow us to complete the level set argument.

\subsection{Organization of the paper} 
The rest of the paper is organized as follows. In Section~\ref{wei-ses}, we recall and collect some important properties of the $A_p$ class of Muckenhoupt weights, and introduce a family of weighted parabolic cylinders and a quasi-distance adapted to the class of weights considered in the paper. In Section~\ref{Fn-spaces}, function spaces and weak solutions of \eqref{main-eqn} are defined, and a compactness lemma for weighted Sobolev spaces is recalled.  The main part of the paper is Section~\ref{bdr-section}. In this section, we first obtain boundary Lipschitz estimates for a simple frozen-coefficient model, then derive boundary energy estimates and boundary approximation estimates for our singular-degenerate equation on flat domains. These are combined with a level-set argument to prove the boundary $L^p$ regularity estimates for $p\geq 2$. Section~\ref{global-section} is devoted to the proof of Theorem~\ref{main-theorem}. Lastly, we end the paper with some admissible weight examples of our results in Section~\ref{Examples}.

\section{Weighted inequalities and weighted parabolic cylinders} \label{wei-ses} 
In this section, we introduce two weighted inequalities and weighted parabolic cylinders. These results will be used in the subsequent sections to establish regularity estimates for solutions to \eqref{main-eqn}.
\subsection{Weights and weighted inequalities}
\begin{definition} \label{A-p-def}  Let $p \in (1, \infty)$, and let $\mu : \mathbb{R}^n \rightarrow [0, \infty)$ be a locally integrable function. We say that $\mu$ belongs to the $A_p$ class of  Muckenhoupt weights if $[\mu]_{A_p} < \infty$, where
\[
[\mu]_{A_p}=\sup_{r>0,\, x_0 \in \mathbb{R}^n}\left(\fint_{B_r(x_0)}\mu(x)\, dx \right)\left(\fint_{B_r(x_0)}\mu(x)^{-\frac{1}{p-1}}\, dx\right)^{p-1}.
\]
\end{definition}
\noindent
Below, we collect and list some important properties of $A_p$ weights that will be used in this paper.  We refer the reader to \cite[Proposition~7.1.5, Corollary~7.2.6, and Proposition~7.2.8]{Grafakos-2} for proof details. 

\begin{proposition}\label{Ap-prop}
Let $p\in (1, \infty)$ and $M_0\geq 1$. 
Suppose $\mu\in A_p$ with $[\mu]_{A_p}\leq M_0$. Then the following properties hold.

\begin{enumerate}[leftmargin=*]
\item[\textup{(i)}]  $[\mu]_{A_p} \geq 1$.

\item[\textup{(ii)}] The weight function $\mu^{-\frac{1}{p-1}}$ is in $A_{p'}$ with $p'=\frac{p}{p-1}\in (1, \infty)$, and it satisfies
\[
[\mu^{-\frac{1}{p-1}}]_{A_{p'}}=[\mu]_{A_p}^{\frac{1}{p-1}}\leq M_0^{\frac{1}{p-1}}.
\]
Therefore, $\mu\in A_2$ if and only if $\mu^{-1}\in A_2$ and $[\mu]_{A_2}=[\mu^{-1}]_{A_2}$. 

\item[\textup{(iii)}] The induced measure $\mu$ is doubling:
\[
\mu(B_{\lambda r}(x)) \leq \lambda^{np}M_0 \mu(B_r(x)) \quad \text{for all}\quad \lambda > 1,\, r>0,\, x\in \mathbb{R}^n.
\]

\item[\textup{(iv)}] There exist constants $\gamma = \gamma(n,p, M_0)>0$ sufficiently small and $N = N(n, p, M_0)\geq 1$ such that $\mu \in A_{p-\gamma}$ and
\[ 
[\mu]_{A_{p-\gamma}} \leq N.
\]

\item[\textup{(v)}] There exist constants $N_2 = N_2(n, p, M_0) > 0$ and $\zeta_0 = \zeta_0(n, p,  M_0) \in (0,1)$ such that
\begin{equation*} \label{reverse holder}
\mu(S)\leq N_2\left(\frac{|S|}{|B|}\right)^{\zeta_0}\mu(B)\quad\text{for every ball }B\subset\R^n\text{ and measurable }S\subset B.
\end{equation*}
\end{enumerate}
\end{proposition}

\begin{remark}\label{A-p-remark} Let  $\beta$  satisfy \eqref{beta-cond}, and let $\bar{\beta}(x) = \beta(x)^{\frac{n_0}{2}}$ for $n_0=\max\{n, 2\}$. It follows from \textup{(ii)} of Proposition~\ref{Ap-prop} that 
\[
\bar{\beta} \in A_{1+\frac{n_0}{2}} \quad \text{with}\quad
[\bar{\beta}]_{A_{1+\frac{n_0}{2}}}\leq M_0^{\frac{n_0}{2}}. 
\]
Consequently, $\bar{\beta}$ also satisfies the properties in Proposition~\ref{Ap-prop}. 
\end{remark}

Note also that
\[
\mathbb{R}^n_+ = \mathbb{R}^{n-1} \times (0, \infty).
\]
The following interpolation lemma will be used in the boundary local energy estimates for the equations considered in the paper. We refer the reader to \cite[Lemma~2.7]{FP1} for the complete proof.

\begin{lemma} \label{embedd-lemma} For every $M_0 \geq 1$, there exist constants $\theta =\theta(n, M_0) \in (0,1)$ sufficiently close to $1$, and $N = N(n, M_0)>0$ such that 
\[
\begin{split}
& \frac{1}{(\beta)_{B_r(x_0)}}\fint_{Q} \beta(x)u(x,t)^2 dx dt \\ 
& \leq N  \left(\fint_{Q} u(x,t)^2 dxdt \right)^{1-\theta} \left[ \left(\fint_{Q} u(x,t)^{2} dxdt \right)^{\theta} + r^{2\theta} \left(\fint_{Q} |\nabla u(x,t)|^{2} dxdt \right)^{\theta} \right],
\end{split}
\]
for every weight $\beta$ satisfying \eqref{beta-cond}, for $u\in L^2(\Gamma, W^{1,2}(B_r(x_0)))$ with some open interval $\Gamma\subset \R$, 
$Q= B_r(x_0) \times \Gamma$ with some $r>0$ and $x_0 \in \R^n$.  The same assertion holds for $u\in L^2(\Gamma, W^{1,2}(B_r^+(x_0)))$ and $Q=B_r^+(x_0)\times\Gamma$ when $x_0\in\overline{\R^n_+}$.
\end{lemma}

Next, we introduce the following inequality that is needed in the paper. 

\begin{lemma}\label{embedd-lemma-2}  For every $M_0 \geq 1$ and $\epsilon \in (0,1)$, there exists $N=N(n,M_0,\epsilon)>0$ such that
\begin{align*}
\fint_{B_r}|u(x)|^2 dx\leq \epsilon\, r^2\fint_{B_r}|\nabla u(x)|^2 dx + \frac{N}{\mu(B_r)}\int_{B_r} |u(x)|^2 \mu(x) dx 
\end{align*}
for every $u \in W^{1,2}(B_r)\cap L^2(B_r, \mu)$ with $\mu\in A_2$ satisfying $[\mu]_{A_2}\leq M_0$, and with $r>0$. 
\end{lemma}

\begin{proof} Without loss of generality, we assume $r =1$. From \textup{(iv)} of Proposition~\ref{Ap-prop}, there exists a sufficiently small constant $\gamma = \gamma(n, M_0) \in (0, 1)$ such that $\mu \in A_{2-\gamma}$, and
\begin{equation*} 
[\mu]_{A_{2-\gamma}} \leq N_0, \quad \text{where} \quad N_0 = N_0(n, M_0) \geq 1.
\end{equation*}
Let us denote $q = \frac{2}{2-\gamma} \in (1,2)$. We write
\[
\fint_{B_1} |u(x)|^{q} dx = \fint_{B_1} \left[|u(x)|^2\mu(x)\right]^{\frac{1}{2-\gamma}} \mu(x)^{-\frac{1}{2-\gamma}}dx,
\]
and apply H\"{o}lder's inequality to obtain
\begin{align} \notag
\left(\fint_{B_1} |u(x)|^{q} dx \right)^{\frac{1}{q}}& \leq \left(\fint_{B_1} |u(x)|^2 \mu(x) dx \right)^{\frac{1}{2}} \left(\fint_{B_1} \mu(x)^{-\frac{1}{1-\gamma}} dx \right)^{\frac{1- \gamma}{2}} \\ \notag
& \leq [\mu]_{A_{2-\gamma}}^{\frac{1}{2}}\left(\fint_{B_1}\mu(x)dx\right)^{-\frac{1}{2}} \left(\fint_{B_1} |u(x)|^2 \mu(x) dx \right)^{\frac{1}{2}} \\ \label{wei-hold-0806}
& \leq \left(\frac{N_0}{\mu(B_1)}\int_{B_1} |u(x)|^2 \mu(x) dx \right)^{\frac{1}{2}}.
\end{align}
On the other hand, by the Gagliardo-Nirenberg inequality, we have
\[
\|u\|_{L^2(B_1)} \leq N \|\nabla u\|_{L^2(B_1)}^{\theta} \|u\|_{L^q(B_1)}^{1-\theta} + N\|u\|_{L^q(B_1)},
\]
where $N = N(n, q)>0$ and $\theta \in (0,1)$ satisfying
\[
\frac{1}{2} = \theta \Big( \frac{1}{2} - \frac{1}{n} \Big) + \frac{1-\theta}{q}.
\]
Then, the assertion of the lemma follows from \eqref{wei-hold-0806}, the last inequality, and Young's inequality. The proof is completed.
\end{proof}

\subsection{Weighted parabolic cylinders, and a quasi-distance function}  \label{cylinder-distance-sec}
Let $M_0 \geq 1$, and let $\beta: \mathbb{R}^n \rightarrow [0, \infty)$ be a weight satisfying \eqref{beta-cond}. We give the definition of weighted parabolic cylinders introduced in \cite{CFP, FP1} that are intrinsically suitable for the class of equations \eqref{main-eqn}. 
\begin{definition}
Given $z_0 = (x_0, t_0) \in \R^{n}\times \R$ and $r>0$, we define the weighted parabolic cylinder $Q_{r,\beta}(z_0)$ of radius $r$ centered at $z_0$, with respect to the weight $\beta$, by
\begin{equation*} 
Q_{r, \beta}(z_0) =B_r(x_0)  \times \Gamma_{\beta, z_0}(r), \quad \text{where}   \quad \Gamma_{\beta, z_0} (r) =(t_0 -r^2\mathcal{E}_{\beta}(x_0, r), t_0],
\end{equation*}
in which $\mathcal{E}_{\beta}(x_0, r)$ is defined by
\begin{equation*} 
\mathcal{E}_{\beta}(x_0, r)= \Big( \fint_{B_{r}(x_0)} \beta(x)^{\frac{n_0}{2}} dx \Big)^{\frac{2}{n_0}},\quad \text{with}\quad n_0=\max\{n, 2\}.
\end{equation*}
\end{definition}
\noindent
By a simple computation, we infer
\begin{equation}\label{heights}
   r^2 \mathcal{E}_{\beta}(x_0, r)=
   \left\{
   \begin{aligned}
   &\tfrac{1}{2}r\bar{\beta}(B_r(x_0)) \quad &&\text{when}\quad n=1, \\
&\sigma_n^{-\frac{2}{n}} \big[\bar{\beta}(B_r(x_0))\big]^{\frac{2}{n}}, \quad &&\text{when}\quad n\geq 2,
   \end{aligned}\right.
\end{equation}
where $\sigma_n$ ($n\geq 2$) denotes the Lebesgue measure of the unit ball $B_1\subset \R^n$, and $\bar{\beta}(x) =\beta(x)^{\frac{n_0}{2}}$, and 
\[
\bar{\beta}(B_r(x_0)) = \int_{B_r(x_0)} \bar{\beta}(x)\, dx.
\]
It follows from \eqref{heights} that the class of weighted cylinders satisfies an inclusion property:
\[
Q_{r_1,\, \beta}(z_0) \subset Q_{r_2,\, \beta}(z_0), \quad \text{for\ } \, 0<r_1 \leq r_2, \  \forall\, z_0 \in \mathbb{R}^{n+1}.
\]
This inclusion only requires the local integrability of $\beta^{n_0/2}$; condition \eqref{beta-cond} also provides the quantitative doubling properties used below. Moreover, from this point of view, the $A_{1+2/n_0}$ class of weights is the class that allows the weighted cylinders, defined in the way of \eqref{cylinder-def}, to satisfy the above inclusion property.

\smallskip
In addition, we denote the upper half ball in $\mathbb{R}^n$ by
\begin{equation}\label{half ball}
B_r^+(x_0) = B_r(x_0) \cap \mathbb{R}^n_+ \quad \text{with} \quad \mathbb{R}^n_+ =  \mathbb{R}^{n-1} \times (0, \infty),
\end{equation}
and the flat boundary portion of $B_r^+(x_0)$ is denoted by
\begin{equation} \label{flat-bdr}
T_r(x_0) =  B_r(x_0) \cap \big(\mathbb{R}^{n-1} \times \{0\}\big).
\end{equation}
The upper half parabolic weighted cylinder of radius $r$ centered at $z_0=(x_0, t_0)$ is denoted by
\begin{equation}\label{half cylinder}
Q_{r, \beta}^+(z_0) = B_r^+(x_0) \times \Gamma_{\beta, z_0}(r).
\end{equation}
Throughout the paper, when $z_0 =(0, 0)$ and $x_0 = 0$, we omit $z_0$ and $x_0$ in the notation. For example, we write
\[
 Q_{r, \beta} = Q_{r, \beta}(0), \quad B_r = B_r(0),\quad T_r = T_r(0), \quad\mathcal{E}_{\beta}(r) = \mathcal{E}_{\beta}(0, r), \quad \Gamma_{\beta}(r) = \Gamma_{\beta, 0} (r).
\]
The same notation applies to $Q_{r, \beta}^+$ and $B_r^+$.

\smallskip
Next, we introduce a quasi-distance adapted to the parabolic weighted cylinders. To begin with, given $z_0=(x_0,t_0)\in\R^{n}\times \R$, we let $h_{x_0}: [0,\infty)\rightarrow [0,\infty)$ be defined by
\begin{equation*}
h_{x_0}(r)=r^2\mathcal{E}_{\beta}(x_0, r).
\end{equation*}
For $r>0$, $\mathcal E_\beta(x_0,r)$ is given in \eqref{t-height}, and we set $h_{x_0}(0)=0$. By \eqref{heights}, positivity of the weight, and \textup{(v)} of Proposition~\ref{Ap-prop}, the function $h_{x_0}$ is continuous and strictly increasing, with $h_{x_0}(r)\to\infty$ as $r\to\infty$. Thus its inverse is defined on $[0,\infty)$. Now, for any $z=(x,t)\in \R^{n} \times \R$ with $t\leq t_0$, we denote
\begin{equation*}
\begin{aligned}
\tilde \rho_{\beta}(z,z_0)
=\inf\big\{r>0: z\in Q_{r,\beta}(z_0)\big\}=\max\big\{|x-x_0|,\ h_{x_0}^{-1}(t_0-t)\big\}.
\end{aligned}
\end{equation*}
Then, for  any $z=(x,t)\in \R^{n}\times \R$, we define
\begin{equation}\label{rho def}
\rho_{\beta}(z, z_0)=
\left\{
\begin{array}{rl}
\tilde \rho_{\beta}(z, z_0) \quad &\text{if}\quad t\leq t_0,\\[4pt]
\tilde \rho_{\beta}(z_0, z)\quad &\text{if}\quad t>t_0.
\end{array}\right.
\end{equation}
Intuitively, $\rho=\rho_{\beta}(z, z_0)$ represents the spatial radius for which $z$ lies on the parabolic boundary $\partial'Q_{\rho,\beta}(z_0)$ if $t \leq t_0$ or $z_0$ lies on the parabolic boundary $\partial'Q_{\rho,\beta}(z)$ when $t >t_0$. 
From the definition in \eqref{rho def}, it is straightforward to see that
\begin{equation}\label{d-r}
|x-x_0|\leq\rho_{\beta}(z, z_0),\quad \text{for all}\quad z=(x,t),\ z_0=(x_0,t_0)\in \R^{n} \times \R.
\end{equation}
Additionally, for any $z_0 = (x_0, t_0) \in \R^{n}\times \R$ and $r>0$, we also have
\begin{equation}\label{cylinder-quasi}
Q_{r,\beta}(z_0) = \{z = (x,t) \in \R^{n}\times \R: \quad \rho_{\beta} (z_0, z) <r \ \text{ and }\  t  \leq t_0\}.
\end{equation}

\smallskip
The following lemma states that the function $\rho_\beta(z, z_0)$ defined in \eqref{rho def} is a quasi-metric on $\R^{n}\times \R$. We refer the reader to \cite[Appendix A]{FP1} for the proof.

\begin{lemma}\label{quasi-metric lemma} Suppose that \eqref{beta-cond} holds. Then, the function $\rho_\beta$ defined in \eqref{rho def}  is a quasi-metric. In particular, there exists a constant $\Lambda=\Lambda(n,M_0)\geq2$ such that
\begin{equation*}\label{quasi-tri}
\rho_{\beta}(z_0,\bar{z})\leq \Lambda \big[\rho_{\beta}(z_0,z_1)+\rho_{\beta}(z_1,\bar{z})\big], \quad \forall z_0,\ z_1,\ \bar{z} \ \in \ \R^{n} \times \R.
\end{equation*}
\end{lemma}

\section{Function Spaces, Weak Solutions, and a Compactness Result} \label{Fn-spaces}
\subsection{Function spaces and weak solutions}\label{scaling-sec}
Let $p \in (1, \infty)$, and $\Omega \subset \mathbb{R}^n$ be an open bounded set. As usual,  $W^{1,p}(\Omega)$ denotes the Sobolev space consisting of functions $u: \Omega \rightarrow \mathbb{R}$ such that $u \in L^p(\Omega)$ and the weak derivatives $u_{x_i} \in L^p(\Omega)$ for all $i =1,2,\ldots, n$, and
\[
\|u\|_{W^{1,p}(\Omega)}^p = \|u\|_{L^p(\Omega)}^p + \|\nabla u\|_{L^p(\Omega)}^p.
\]
For real numbers $S<T$, we denote by $L^p((S,T),W^{1,p}(\Omega))$ the Bochner space consisting of all measurable functions $u: (S, T) \rightarrow W^{1, p}(\Omega)$ such that 
\[
\|u\|_{L^p((S, T), W^{1,p}(\Omega))} = \left( \int_S^T \|u(\cdot, t)\|_{W^{1,p}(\Omega)}^p\, dt \right)^{1/p} < \infty.
\]
We also denote by $W^{1,p}(\Omega)^*$ the dual space of $W^{1,p}(\Omega)$. It is well-known that every functional $\ell\in(W^{1,p'}(\Omega))^*$ has representations as
\[
\langle\ell,\varphi\rangle=\int_\Omega\left(f(x)\varphi(x)+F(x)\cdot\nabla\varphi(x)\right)\,dx,
\qquad \forall \ \varphi \ \in \ W^{1,p'}(\Omega),
\]
for some $f\in L^p(\Omega)$ and $F\in L^p(\Omega)^n$, and its norm is given by
\[
\|\ell\|_{W^{1,p'}(\Omega)^*}=\inf_{(f, F)}\left(\|f\|_{L^p(\Omega)}^p+\|F\|_{L^p(\Omega)}^p\right)^{1/p},
\]
where the infimum is taken over all representations of $\ell$ as above. A similar definition is also formulated for the Bochner space $L^p((S, T), (W^{1,q}(\Omega))^{*})$ with $p, q \in (1, \infty)$.  We now define the weighted Sobolev spaces used in the paper.

\begin{definition} \label{W-space} Let $p \in (1, \infty)$, $\Omega \subset \mathbb{R}^n$ be a non-empty open bounded set, and $S<T$. For a given locally integrable function $\beta: \Omega\rightarrow [0, \infty)$, a function $u\in L^p((S,T), W^{1,p}(\Omega))$ is said to belong to the parabolic weighted Sobolev space $\W^{1,p}(\Omega\times (S,T), \beta)$ if $\beta u\in L^1_{\textup{loc}}(\Omega\times (S,T))$ and 
\[
\beta u_t \in L^p((S, T), W^{1, p'}(\Omega)^{*}),
\] 
where $p'=p/(p-1)\in (1, \infty)$. The space $\W^{1,p}(\Omega\times (S,T), \beta)$ is equipped with the norm 
\begin{equation}\label{Norm-W1p}
\|u\|_{\W^{1,p}(\Omega\times (S,T), \beta)} = \|u\|_{L^p((S,T), W^{1,p}(\Omega))} + \| \beta u_t\|_{L^p((S,T),\, W^{1,p'}(\Omega)^{*})}.
\end{equation}
\end{definition} 
\noindent
Now, let $\beta$ be a weight satisfying \eqref{beta-cond}, and $\bA$ be as in \eqref{ellip-cond}. We introduce the definitions of weak solutions that are needed for the theory of interior estimates for \eqref{main-eqn}.

\begin{definition} \label{interior-weak-def} Let $p \in (1, \infty)$, $\Omega\subset \R^n$ be an open bounded domain with boundary $\partial\Omega \in C^1$, and $Q = \Omega \times (s, \tau)$ for two given numbers $s<\tau$. 
\begin{enumerate}[leftmargin=*]
\item[\textup{(i)}] 
Given $F\in L^p(Q)^n$, a function $u\in L^p_{\mathrm{loc}}((s,\tau),W^{1,p}_{\mathrm{loc}}(\Omega))$ with $\beta u\in L^1_{\mathrm{loc}}(Q)$ is said to be a weak solution of
\begin{equation} \label{scaling-eqn}
\beta(x) u_t - \textup{div}(\bA \nabla u) = \textup{div}(F)\quad \text{in} \quad Q,
\end{equation}
if  for every $\varphi \in C_0^\infty(Q)$, we have
\begin{align*}
-\int_{Q} \beta(x) u(x,t) \varphi_t(x,t) dx dt &+ \int_{Q} \wei{\bA(x,t) \nabla u(x,t), \nabla \varphi(x,t)} dx dt \\&= -\int_{Q} \wei{F(x,t), \nabla \varphi(x,t)} dx dt.
\end{align*}

\item[\textup{(ii)}] Given $F\in L^p(Q)^n$ and $g\in L^p(\partial \Omega\times (s,\tau))$, a function $u \in \W^{1,p}(Q, \beta)$ is said to be a weak solution of
\begin{equation}\label{soln-def}
\left\{
\begin{aligned}
\beta(x) u_t - \textup{div}(\bA \nabla u) &= \textup{div}(F)\quad &&\text{in} \quad Q,\\
(\bA\nabla u+F)\cdot \vec{\mathbf{n}} &= g \quad &&\text{on} \quad \partial \Omega \times (s,\tau),
\end{aligned}\right.
\end{equation}
if for every $\varphi \in C^\infty(\overline{Q})$ vanishing near $t=s$ and $t=\tau$, we have
\begin{align*}
-\int_{Q} \beta(x) u(x,t) \varphi_t(x,t) dx dt &+ \int_{Q} \wei{\bA(x,t) \nabla u(x,t), \nabla \varphi(x,t)} dx dt\\& = -\int_{Q} \wei{F(x,t), \nabla \varphi(x,t)} dx dt+\int_{\partial\Omega\times(s,\tau)}g(x,t)\varphi(x,t)\,dSdt.
\end{align*}
\end{enumerate}
\end{definition}

\smallskip
We next introduce needed function spaces and the notion of weak solutions to study local boundary estimates for \eqref{main-eqn}. Let $r>0$ and $x_0\in\overline{\mathbb R^n_+}$ be such that $\partial B_r(x_0)\cap\mathbb R^n_+$ is non-empty. For any $p\in (1, \infty)$, we first define a subspace of $W^{1,p}(B_r^+(x_0))$:
\begin{equation}\label{space-Vp}
\mV_p(B_r^+(x_0))=\left\{\varphi\in W^{1,p}(B_r^+(x_0)): \  \left.\varphi\right|_{\partial B_r(x_0)\cap\mathbb R^n_+}=0
\quad\text{in the sense of traces} \right\}.
\end{equation}
This is a closed subspace, equipped with the norm
\[
\|\varphi\|_{\mV_p(B_r^+(x_0))}=\|\varphi\|_{W^{1,p}(B_r^+(x_0))}.
\]
Note also that functions in $\mV_p(B_r^+(x_0))$ have no trace condition on the flat part $T_r(x_0)$. 

\begin{definition}\label{space-hW1p}
Let $S<T$, $p'=p/(p-1)$, and let $\beta$ be a nonnegative
locally integrable function on $B_r^+(x_0)$.
We define $\hW^{1,p}(B_r^+(x_0)\times(S,T),\beta)$ as the space consisting of functions
\[
u\in L^p((S,T), W^{1,p}(B_r^+(x_0)))
\]
such that
\[
\beta u\in L^1_{\textup{loc}}(B_r^+(x_0)\times(S, T)) \quad \text{and}\quad \beta u_t\in L^p((S,T), \mV_{p'}(B_r^+(x_0))^*),
\]
where $\mV_{p'}(B_r^+(x_0))^*$ is the dual space of $\mV_{p'}(B_r^+(x_0))$. The space $\hW^{1,p}(B_r^+(x_0)\times(S,T),\beta)$ is equipped with the norm
\begin{align*}
\|u\|_{\hW^{1,p}(B_r^+(x_0)\times(S,T),\beta)}
=\|u\|_{L^p((S,T), W^{1,p}(B_r^+(x_0)))}+\|\beta u_t\|_{L^p((S,T), \mV_{p'}(B_r^+(x_0))^*)}.
\end{align*}
When $\beta$ is a positive constant, we simply write
\[
\hat{\W}^{1,p}(B_r^+(x_0)\times(S,T))=\hat{\W}^{1,p}(B_r^+(x_0)\times(S,T),\beta).
\]
\end{definition}

\begin{definition} \label{boundary-weak-def} Let $p \in (1, \infty)$, $r>0$ and $z_0 = (x_0, t_0) \in \overline{\mathbb R^n_+}\times\mathbb R$. We write $Q^+ = Q_{r, \beta}^+(z_0)$ and suppose that $T_r(x_0) \not= \emptyset$. Given $F\in L^p(Q^+)^n$ and $g\in L^p(T_r(x_0) \times \Gamma_{\beta, z_0} (r))$, a function $u \in \hW^{1,p}(Q^+,\beta)$ is said to be a weak solution of
\begin{equation}\label{bdr-scaling-eqn}
   \left\{
\begin{aligned}
\beta(x) u_t - \textup{div}(\bA\nabla u) &= \textup{div}(F) \quad &&\text{in} \quad Q^+,\\
(\bA\nabla u+F)\cdot \vec{\mathbf{n}} &= g \quad &&\text{on} \quad T_r(x_0) \times \Gamma_{\beta, z_0} (r),
\end{aligned}\right.
\end{equation}
if 
\begin{align*}
& -\int_{Q^+} \beta(x) u(x,t) \varphi_t(x,t) dx dt + \int_{Q^+} \wei{\bA(x,t) \nabla u(x,t), \nabla \varphi(x,t)} dx dt \\
& = -\int_{Q^+} \wei{F(x,t), \nabla \varphi(x,t)} dx dt + \int_{T_r(x_0) \times \Gamma_{\beta, z_0} (r)} g(x', t) \varphi(x', 0, t) dx' dt,
\end{align*}
for every  $\varphi\in C_0^\infty(B_r(x_0)\times(t_0-r^2\mathcal E_\beta(x_0,r),t_0))$.
\end{definition}

Based on the definitions above, we discuss some scaling properties of the class of equations \eqref{main-eqn} here. Besides the standard heat dilation $(x,t) \mapsto (rx, r^2t)$, let us introduce two special dilation properties that will be used frequently throughout the paper. Let us consider the equation \eqref{scaling-eqn}. After taking a translation, we may assume without loss of generality that $z_0 = (0,0)$. Recall that $\mathcal{E}_\beta (r) = \mathcal{E}_{\beta}(0, r)$, as defined in \eqref{t-height}.  For $(x,t) \in Q_1=B_1\times (-1,0)$, let us now define the following:
\begin{equation} \label{scale-1}
\left\{
\begin{aligned}
\tilde{u}(x,t) &= u(rx, r^2\mathcal{E}_\beta(r) t),
\quad &&\tilde{\beta}(x) = [\mathcal{E}_\beta(r)]^{-1} \beta(rx),\\[4pt]
\tilde{\bA} (x,t)&= \bA(rx, r^2\mathcal{E}_\beta(r) t), \quad
&&\tilde{F}(x,t) =rF(rx, r^2\mathcal{E}_\beta(r) t),\\
\tilde{g}(x, t) &= r g(rx, r^2\mathcal{E}_\beta(r) t).
\end{aligned} \right.
\end{equation}
In addition, the class of equations \eqref{main-eqn} is specially invariant under the following time only dilation:
\begin{equation} \label{scale-2}
\left\{
\begin{aligned}
\hat{u}(x,t) &= u(x,\mathcal{E}_{\beta}(r) t),
\quad &&\hat{\beta}(x) = [\mathcal{E}_\beta(r)]^{-1} \beta(x),\\[4pt]
\hat{\bA} (x,t)&= \bA(x, \mathcal{E}_\beta(r) t), \quad
&&\hat{F}(x,t) = F(x,\mathcal{E}_{\beta}(r) t),\quad \hat{g}(x, t) = g(x,\mathcal{E}_{\beta}(r) t).
\end{aligned}  \right.
\end{equation}
It can be verified directly from \eqref{t-height} that
\[
\mathcal{E}_{\tilde \beta}(1)=1 \quad \text{and} \quad \mathcal{E}_{\hat{\beta}}(r) =1 \quad \text{so that}\quad Q_{1,\tilde \beta} = Q_1 \quad \text{and}  \quad Q_{r,\hat{\beta}} = Q_r.
\]
\noindent
The following lemma will be used frequently throughout the paper. Its proof follows from a direct calculation, and is therefore omitted.
\begin{lemma} \label{scaling-lemma} Let $F\in L^p(Q_{r,\beta}^+)^n$ and $g\in L^p(T_r\times\Gamma_\beta(r))$ for some $p \in (1, \infty)$, and suppose that $u \in \hW^{1,p}(Q_{r,\beta}^+, \beta)$ is a weak solution of \eqref{bdr-scaling-eqn} in the sense of Definition~\ref{boundary-weak-def}. Then,
\begin{enumerate}[leftmargin=*]
\item[\textup{(i)}]    
$\tilde u \in \hW^{1,p}(Q_{1}^+, \tilde{\beta})$ is a weak solution of
\begin{equation*}
\left\{
\begin{aligned}
\tilde{\beta}(x) \tilde{u}_t - \textup{div}(\tilde{\bA} \nabla \tilde{u}) &= \textup{div}(\tilde{F}) \quad &&\text{in} \quad Q_1^+,\\
(\tilde{\bA}\nabla \tilde{u}+\tilde{F})\cdot \vec{\mathbf{n}} &= \tilde{g} \quad &&\text{on} \quad T_1 \times (-1,0),
   \end{aligned}\right.
\end{equation*}
where $\tilde{u}$, $\tilde{\beta}$, $\tilde{\bA}$, $\tilde{F}$ and $\tilde g$ are defined in \eqref{scale-1}. 

\item[\textup{(ii)}]$\hat{u} \in  \hW^{1,p}(Q_{r}^+, \hat{\beta})$ is a weak solution of
\[
\left\{
\begin{aligned}
\hat{\beta}(x) \hat{u}_t - \textup{div}(\hat{\bA} \nabla \hat{u}) &= \textup{div}(\hat{F}) \quad &&\text{in} \quad Q_{r}^+,\\
(\hat{\bA}\nabla \hat{u}+\hat{F})\cdot \vec{\mathbf{n}} &= \hat{g} \quad &&\text{on} \quad T_r \times (-r^2, 0),
\end{aligned}\right.
\]
where $\hat u$, $\hat{\beta}$, $\hat{\bA}$, $\hat{F}$ and $\hat{g}$ are defined in \eqref{scale-2}.
\end{enumerate} 
\end{lemma}

\subsection{A compactness theorem}
In this subsection, we introduce a compactness result, Theorem~\ref{compactness-lemma}, for weak solutions of \eqref{soln-def} in $Q_{r,\beta}^+$. To this end, we state and prove a lemma that is needed in later sections.

\begin{lemma}\label{beta-seq}
Let $M_0\geq 1$, and denote $n_0 = \max\{n, 2\}$. There exists a constant $q_0 = q_0(n, M_0)\in (1, 2)$ such that the following assertions hold. Suppose that $\{\beta_k\}$ is a sequence of non-negative locally integrable functions defined in $\mathbb{R}^n$ satisfying \eqref{beta-cond} and
\begin{equation}  \label{compact-beta-ass}
\Big( \fint_{B_r(x_0)} \beta_k(x)^{\frac{n_0}{2}}\, dx \Big)^{\frac{2}{n_0}}=1, \quad \forall \ k \in \mathbb{N}
\end{equation}
for some given $r>0$ and $x_0\in\R^n$, and also that 
\begin{equation}\label{osc-assumption}
\lim_{k \rightarrow \infty} \frac{1}{\beta_k(B_r(x_0))}\int_{B_r(x_0)} |\beta_k(x) - (\beta_k)_{B_r(x_0)}|^2 \beta_k(x)^{-1}\, dx =0.
\end{equation}   
Then, there exists a constant $\beta_0\in [\frac{1}{M_0},1]$ such that, up to a subsequence, we have
\begin{equation}\label{beta-lim}
\lim_{k\rightarrow \infty}(\beta_k)_{B_r(x_0)}=\beta_0, \quad \text{and}\quad 
\beta_k(x) \rightarrow \beta_0 \quad \text{in} \quad L^{q_0}(B_r(x_0)) \quad \text{as} \quad k \rightarrow \infty.
\end{equation}
\end{lemma}

\begin{proof} By translation $x \mapsto x - x_0$, we can assume without loss of generality that $x_0 =0$. Note from \eqref{beta-cond}, \eqref{compact-beta-ass}, and the definition of $A_{1+\frac{2}{n_0}}$ weights, we have 
\[
(\frac{1}{\beta_k})_{B_r} = \fint_{B_r} \beta_k(x)^{-1}\, dx  = \Big( \fint_{B_r} \beta_k(x)^{-1}\, dx \Big)\Big( \fint_{B_r} \beta_k(x)^{\frac{n_0}{2}}\, dx \Big)^{\frac{2}{n_0}} \leq M_0.
\]
Note also that we have $1 \leq (\frac{1}{\beta_k})_{B_r} (\beta_k)_{B_r}$ due to H\"{o}lder's inequality. Then, by using H\"{o}lder's inequality again when $n \geq 3$, we have
\begin{equation} \label{M-0-evar-0805}
\frac{1}{M_0}\leq [(\frac{1}{\beta_k})_{B_r}]^{-1}\leq(\beta_k)_{B_r}\leq \Big( \fint_{B_r} \beta_k(x)^{\frac{n_0}{2}}\, dx \Big)^{\frac{2}{n_0}}=1, \quad \forall \ k\in\N.
\end{equation}
This means that the sequence of numbers $\{(\beta_k)_{B_r}\}$ is bounded. Therefore, there exists a constant $\beta_0\in [\frac{1}{M_0},1]$ such that, up to a subsequence, we have
\begin{equation}\label{beta-lim-1}
\lim_{k\rightarrow \infty}(\beta_k)_{B_r}=\beta_0. 
\end{equation}
This proves the first assertion of \eqref{beta-lim}.
Besides, we also note from \eqref{M-0-evar-0805} that
\begin{equation}\label{ball-measures}
\beta_k(B_r)\leq |B_r| \leq \beta_k^{-1}(B_r) \quad \text{for all}\quad k\in \mathbb{N}, \quad \text{where} \quad \beta_k^{-1}(B_r) = \int_{B_r} \beta_k(x)^{-1}\, dx.
\end{equation}

\smallskip
Due to \eqref{beta-cond} and by \textup{(iv)} of Proposition~\ref{Ap-prop}, there exist a sufficiently small constant $\gamma = \gamma(n, M_0) \in (0, \frac{2}{n_0})$ and a constant  $N_0 = N_0(n, M_0)\geq 1$ such that 
\begin{equation}\label{mu-p-gamma-1}
\frac{1}{\beta_k} \in A_{1+\frac{2}{n_0}-\gamma} \quad \text{with}\quad \big[\frac{1}{\beta_k}\big]_{A_{1+\frac{2}{n_0}-\gamma}} \leq N_0, \quad \text{for all}\quad k\in \mathbb{N}.
\end{equation}
Based on this $\gamma$, we denote $q_0=2/(1+\frac{2}{n_0}-\gamma)\in (1, 2)$, and write
\[
\fint_{B_r} |\beta_k(x) - (\beta_k)_{B_r}|^{q_0} dx = \fint_{B_r} \left[|\beta_k(x) - (\beta_k)_{B_r}|^2\beta_k(x)^{-1}\right]^{\frac{q_0}{2}} \beta_k(x)^{\frac{q_0}{2}}dx.
\]
Then, applying H\"{o}lder's inequality together with \eqref{mu-p-gamma-1}, we get
\begin{equation}\label{q_0 approx}
\begin{aligned}
&\left(\fint_{B_r}  |\beta_k(x) - (\beta_k)_{B_r}|^{q_0} dx \right)^{1/q_0}\\ 
& \leq \left(\fint_{B_r} |\beta_k(x) - (\beta_k)_{B_r}|^2 \beta_k(x)^{-1} dx \right)^{1/2} \left(\fint_{B_r} \beta_k(x)^{\frac{q_0}{2-q_0}} dx \right)^{(2-q_0)/(2q_0)} \\
& \leq \big[\frac{1}{\beta_k}\big]_{A_{1+\frac{2}{n_0}-\gamma}}^{\frac{1}{2}}\left(\fint_{B_r}\beta_k(x)^{-1}dx\right)^{-1/2} \left(\fint_{B_r} |\beta_k(x) - (\beta_k)_{B_r}|^2 \beta_k(x)^{-1} dx \right)^{1/2} \\
& \leq N_0^{\frac{1}{2}}\left(\frac{1}{\beta_k(B_r)}\int_{B_r} |\beta_k(x) - (\beta_k)_{B_r}|^2 \beta_k(x)^{-1}dx \right)^{1/2}\rightarrow 0,\quad \text{as}\quad k\rightarrow \infty,
\end{aligned}
\end{equation}
where in the last inequality we also used \eqref{ball-measures}. By the triangle inequality, we note that
\begin{align*}
\|\beta_k(x)-\beta_0\|_{L^{q_0}(B_r)}
\leq \|\beta_k(x)-(\beta_k)_{B_r}\|_{L^{q_0}(B_r)}
+\|(\beta_k)_{B_r}-\beta_0\|_{L^{q_0}(B_r)}.
\end{align*}
Thus, letting $k\rightarrow \infty$, the second assertion of \eqref{beta-lim}  follows from \eqref{beta-lim-1} and \eqref{q_0 approx}.
\end{proof}

Now, we state the following Aubin-Lions type compactness result introduced in \cite[Proposition~3.4]{FP1} in which the functional spaces vary along the sequences of functions. Here $Q_r^+=B_r^+\times(-r^2,0]$, $2^* = \frac{2n}{n-2}$ when $n \geq 3$ and $2^* = \infty$ when $n=1, 2$. The proof of the theorem can be found in \cite[Proposition~3.4]{FP1}.
\begin{theorem}\label{compactness-lemma} For every $M_0 \geq 1$, there exists a constant $l_0 = l_0(n, M_0)\in (1,2)$ sufficiently close to $1$ such that the following assertions hold. Suppose that $\{\beta_k\}$ is a sequence of non-negative locally integrable functions defined on $\mathbb{R}^n$ each satisfying \eqref{beta-cond} with the same $M_0$, as well as \eqref{compact-beta-ass} and \eqref{osc-assumption} on $B_r$, for some $r>0$. Suppose also that $\{u_k\}$ is a sequence of functions satisfying $u_k\in\hW^{1,2}(Q_r^+,\beta_k)$ for every $k\in\mathbb N$, and there is a constant $M>0$ such that
\begin{equation} \label{u-k-bound-0402}
\|u_k\|_{\hW^{1,2}(Q_r^+,\beta_k)} \leq M, \quad \forall \ k \in \mathbb{N}.
\end{equation}
Then there exists a function $u\in\hW^{1,2}(Q_r^+)$ such that, up to a subsequence, we have
\begin{equation} \label{l-0-convergence}
\left\{
\begin{array}{ll}
\beta_k(u_k)_t\rightharpoonup\beta_0u_t &\quad\text{in}\quad L^2((-r^2,0), \mV_2(B_r^+)^{*}),\\[4pt]
u_k\rightarrow u &\quad\text{in}\quad L^2((-r^2,0),L^{l_0}(B_r^+)),\ \text{and}\\[4pt]
\nabla u_k\rightharpoonup\nabla u &\quad\text{in}\quad L^2(Q_r^+),\qquad\text{as}\quad k\rightarrow\infty,
\end{array}\right.
\end{equation}
where $\beta_0\in [\frac{1}{M_0}, 1]$ is given in Lemma \ref{beta-seq}, and $\mV_2(B_r^+)^{*}$ is the dual space of $\mV_2(B_r^+)$, which is defined in \eqref{space-Vp}. Moreover, for each $l \in [1, 2^*)$, we also have
\begin{equation} \label{L-q-convergence}
u_k\rightarrow u\quad\text{in}\quad L^2((-r^2,0),L^l(B_r^+)),\quad\text{as}\quad k\rightarrow\infty.
\end{equation}
\end{theorem}

\section{Boundary Regularity Estimates on Flat Domains} \label{bdr-section}
The goal of this section is to prove Theorem~\ref{inter-theorem-bdry}. To this end, throughout this section, we assume that $\beta$ is a weight satisfying \eqref{beta-cond} and $\bA$ is a matrix satisfying \eqref{ellip-cond}. In addition, we denote
\begin{equation} \label{Theta-A-def+}
\begin{split}
& \Theta_{\bA, r, x_0} (x,t) = |\bA(x,t) - (\bA)_{B_{r}^+(x_0)}(t)|, \quad 
\Theta_{\bA, r} (x,t)= |\bA(x,t) - (\bA)_{B_{r}^+}(t)|, \\
& \Theta_{\beta, r, x_0} (x) = |\beta(x) - (\beta)_{B_{r}(x_0)}|\beta(x)^{-\frac{1}{2}}, \quad \text{and}\quad 
\Theta_{\beta, r} (x)= |\beta(x) - (\beta)_{B_{r}}| \beta(x)^{-\frac{1}{2}}.
\end{split}
\end{equation}
As such, recalling \eqref{BMO-defi}, we note that
\begin{equation} \label{BMO-A}
[[\bA]]_{\textup{BMO}}(Q_{1, \beta}^+, 1) = \sup_{\substack{z_0=(x_0,t_0) \in Q_{1,\beta}^+\\ r \in (0,1)} }\left(\fint_{Q_{r, \beta}^+(z_0)} \Theta_{\bA, r, x_0} (x,t)^2 dxdt \right)^{1/2},
\end{equation}
and
\begin{equation} \label{BMO-beta}
[[\beta]]_{\textup{BMO}}(B_{1}, 1) = \sup_{x_0 \in B_{1}, r \in (0,1)}\left(\frac{1}{\beta(B_r(x_0))}\int_{B_r(x_0)} \Theta_{\beta, r, x_0}(x)^2  dx\right)^{1/2}.
\end{equation}
\noindent
The proof of  Theorem \ref{inter-theorem-bdry} requires some preparation, and the section is organized as follows. In Subsection~\ref{Lips-est-sec}, we prove boundary Lipschitz estimates for solutions to the equation with frozen coefficients \eqref{v-eqn-bdr}. Local boundary energy estimates and boundary gradient approximation estimates for weak solutions to \eqref{eqn-bdry} are established in Subsection~\ref{bdry energy section} and Subsection~\ref{bdry approx section}, respectively. 
\subsection{Boundary Lipschitz estimates for simple equations}\label{Lips-est-sec}

Let $\bar{\bf{A}}(t)$ be an $n\times n$ matrix function defined on some open interval to be specified. We assume that $\bar{\bf{A}}(t)$ satisfies the ellipticity and boundedness conditions: there exists a constant $\nu \in (0,1)$ such that
\begin{equation} \label{ellip-constant}
\nu |\xi|^2 \leq \wei{\bar{\bA}(t) \xi,\, \xi} \quad \text{and} \quad |\bar{\bA}(t)| \leq \nu^{-1}, \qquad \forall \ \xi \in \mathbb{R}^n.
\end{equation}
In this subsection, we prove the following boundary Lipschitz estimates for the parabolic equation \eqref{main-eqn} with frozen coefficients.

\begin{proposition} \label{Bdry-Lipschitz-constant} For every $\nu \in (0,1)$ and $M_0 \geq 1$, there is a constant $N = N(n, \nu, M_0)>0$ such that the following assertion holds. Let $r>0$, let $\beta(x)$ be a weight satisfying \eqref{beta-cond}, and let $\bar{\bA}(t)$ satisfy \eqref{ellip-constant} on $\Gamma_{\beta}(r)$. Suppose that $v \in \hat{\W}^{1,2}(Q_{r, \beta}^+)$ is a weak solution of
\begin{equation}\label{v-eqn-bdr}
\left\{
\begin{aligned}
(\beta)_{B_r}v_t-\textup{div}(\bar{\bA}(t)\nabla v)&=0 &&\text{in}\quad Q_{r,\beta}^+,\\
(\bar{\bA}(t)\nabla v)\cdot\vec{\mathbf n}&=0 &&\text{on}\quad T_r\times\Gamma_\beta(r),
\end{aligned}\right.
\end{equation}
in the sense of Definition~\ref{boundary-weak-def}. Then
\[
\|\nabla v\|_{L^\infty(Q_{r/2,\beta}^+)}
\leq N\Big(\fint_{Q_{r,\beta}^+}|\nabla v(x,t)|^2dxdt\Big)^{1/2}.
\]
\end{proposition}
\begin{proof} We adapt the argument in \cite[Lemma 5.4]{Dong-Zhang} to our setting, and provide the details here for completeness. As in \eqref{scale-1}, after a dilation, it suffices to prove the proposition for $r=1$ and $\mathcal E_\beta(1)=1$. In this setting, $M_0^{-1}\leq(\beta)_{B_1}\leq1$, and the equation of $v$ reduces to
\[
\left\{
\begin{aligned}
(\beta)_{B_1}v_t-\textup{div}(\bar{\bA}(t)\nabla v)&=0
&&\text{in}\quad Q_1^+,\\
(\bar{\bA}(t)\nabla v)\cdot\vec{\mathbf n}&=0
&&\text{on}\quad T_1\times(-1,0).
\end{aligned}\right.
\]
Thus, this equation is uniformly parabolic. By extending $\bar{\bA}(t)$ to be $\nu I_n$ on $\R\setminus (-1,0)$, then using the mollifier and passing to the limit, we may assume without loss of generality that $\bar{\bA}\in C^\infty((-1,0))$.  Due to this, it follows from local boundary Schauder regularity theory (see \cite[Theorem 4.23]{Lieberman}, for example) that $v \in C^\infty(\overline{Q}_{\rho_0, \beta}^+)$ for all $\rho_0\in(0, 1)$. Then, it remains to prove the estimates when $v \in C^\infty(\overline{Q}_{\rho_0, \beta}^+)$.

\smallskip
Note that for each $i=1, 2, \ldots,n-1$, the function $D_i v$ satisfies the same equation as $v$ with the homogeneous conormal boundary condition on $T_1\times(-1,0)$. Hence, for every $\rho\in(0, \rho_0)$, applying boundary Moser iteration, together with a finite covering, gives
\begin{equation}\label{tangential-derivative}
\sum_{i=1}^{n-1}\|D_i v\|_{L^\infty(Q_{(1+\rho)/2, \beta}^+)}
\leq N(n,\nu,M_0,\rho)
\Big(\fint_{Q_1^+}|\nabla v(x,t)|^2dxdt\Big)^{1/2}.
\end{equation}

Next, we write $\bar{\bA}(t)=(\bar a_{ij}(t))_{i,j=1}^n$ and define
\[
b(x',t) = -\frac{1}{\bar a_{nn}(t)}
\sum_{i=1}^{n-1}\bar a_{ni}(t)D_i v(x', 0, t),
\quad  x' \ \in \ \R^{n-1}, \ |x'| <(1+\rho)/2,  \ t \ \in \ (-1,0).
\]
Now, we let $z=D_nv$. It then follows  that $z$ solves
\[
\left\{
\begin{aligned}
(\beta)_{B_1}z_t-\textup{div}(\bar{\bA}(t)\nabla z)&=0 \quad &&\text{in}\quad Q_{(1+\rho)/2, \beta}^+,\\
z(x', t)&=b(x', t) \quad &&\text{on}\quad T_{(1+\rho)/2}\times\Gamma_\beta((1+\rho)/2),
\end{aligned}\right.
\]
where the flat boundary datum $b$ is bounded due to \eqref{tangential-derivative} and \eqref{ellip-constant}. Indeed, 
\[
\|b\|_{L^{\infty}(T_{(1+\rho)/2}\times\Gamma_\beta((1+\rho)/2))}\leq N_0 \Big(\fint_{Q_1^+}|\nabla v(x,t)|^2dxdt\Big)^{1/2}
\]
for some constant $N_0=N_0(n,\nu,M_0,\rho)>0$. Thus, we can apply boundary Moser iteration to $z$ to infer that
\[
\|z\|_{L^\infty(Q_{\rho, \beta}^+)}\leq N(n,\nu,M_0,\rho)
\left[ N_0 \Big(\fint_{Q_1^+}|\nabla v(x,t)|^2dxdt\Big)^{1/2}+\Big(\fint_{Q_1^+}|z(x,t)|^2dxdt\Big)^{1/2}\right].
\]
Combining the last estimate with \eqref{tangential-derivative}, we obtain
\begin{equation}\label{Lip-k-est-0806}
\|\nabla v\|_{L^\infty(Q_{\rho, \beta}^+)}
\leq N(n,\nu,M_0,\rho)\Big(\fint_{Q_1^+}|\nabla v(x,t)|^2dxdt\Big)^{1/2},
\qquad 0<\rho<\rho_0.
\end{equation}
By the doubling property of $\beta$, we can choose $\rho_0=\rho_0(n,M_0)\in(\frac{1}{2},1)$, close to $1$, such that $Q_{1/2,\beta}^+\subset Q_{(1+\rho_0)/2, \beta}^+$. This and the last estimate with $\rho=\frac{1}{2}$ give the desired estimate. Therefore, the proposition is proved.
\end{proof}

\subsection{Boundary energy estimates}\label{bdry energy section} 
We study the following parabolic problem
\begin{equation} \label{eqn-bdry}
\left\{
\begin{aligned}
\beta(x) u_t - \textup{div}(\bA(x,t) \nabla u)  &= \textup{div}(F(x,t))\quad &&\text{in}\quad Q^+_{r,\beta}(z_0),\\
(\bA(x,t)\nabla u+F)\cdot \vec{\mathbf{n}} &=g(x,t) \quad &&\text{on} \quad T_{r}(x_0) \times \Gamma_{\beta, z_0}(r),
\end{aligned}\right.
\end{equation}
Here $z_0=(x_0,t_0)$ with $x_0\in\overline{\mathbb R^n_+}$, $F\in L^2(Q^+_{r,\beta}(z_0))^n$ and $g\in L^2(T_r(x_0)\times\Gamma_{\beta,z_0}(r))$ are given, and $T_r(x_0)\ne\emptyset$. We note that the weak solutions to \eqref{eqn-bdry} are defined in Definition \ref{boundary-weak-def}. We begin with the following lemma on Caccioppoli type estimates for weak solutions to \eqref{eqn-bdry}.

\begin{lemma} \label{Caccio-2}
For every $\nu\in(0,1)$ and $M_0\geq1$, there exists a constant $N=N(n,\nu,M_0)>0$ such that the following holds. Suppose that $u\in\hat{\W}^{1,2}(Q^+_{r,\beta}(z_0),\beta)$ is a weak solution of \eqref{eqn-bdry}, where $F\in L^2(Q^+_{r,\beta}(z_0))^n$ and $g\in L^2(T_r(x_0)\times\Gamma_{\beta,z_0}(r))$. Then,
\begin{align*}
&\operatorname*{ess\,sup}_{t\in\Gamma_{\beta,z_0}(r)}
\frac{1}{r^2\mathcal E_\beta(x_0,r)}
\fint_{B_r^+(x_0)}\beta(x)u(x,t)^2\varphi(x,t)^2\,dx
+\fint_{Q^+_{r,\beta}(z_0)}
|\nabla(u\varphi)(x,t)|^2\,dxdt\\
&\leq N\left\{
\fint_{Q^+_{r,\beta}(z_0)}|F(x,t)|^2\,dxdt
+\fint_{T_r(x_0)\times\Gamma_{\beta,z_0}(r)}
g(x',t)^2\,dx'dt
\right.\\
&\left.\qquad\quad
+\fint_{Q^+_{r,\beta}(z_0)}
u(x,t)^2\left(1+\frac{1}{r^2}+|\nabla\varphi(x,t)|^2
+\beta(x)|\varphi_t(x,t)|\right)\,dxdt
\right\}
\end{align*}
for every $\varphi\in C^\infty(\overline{Q_{r,\beta}(z_0)})$ satisfying $0\leq\varphi\leq1$ and vanishing near the parabolic boundary $\partial'Q_{r,\beta}(z_0)$.
\end{lemma}

\begin{proof} 
Let $\varphi$ be as in the statement. By using Steklov's average if needed (see \cite[p.~18]{DiB}, for example), we can test the equation \eqref{eqn-bdry} with $u\varphi^2$ to obtain   
\begin{equation}\label{energy-1}
\begin{aligned}
& \frac{1}{2}\frac{d}{dt}\fint_{B_r^+(x_0)} \beta(x)u(x,t)^2 \varphi(x,t)^2  dx + \fint_{B_r^+(x_0)} \wei{\bA \nabla u, \nabla u}(x,t) \varphi(x,t)^2 dx \\
& = -\fint_{B_r^+(x_0)}\Big(2 \wei{\bA \nabla u, \nabla \varphi}u\varphi  + \varphi^2 \wei{F, \nabla u} + 2 u\varphi \wei{F, \nabla \varphi} \Big)(x,t)dx \\
&\quad  + \fint_{B_r^+(x_0)} \beta(x)\big(u^2 \varphi\varphi_t\big)(x,t) dx+ \frac{1}{|B_r^+(x_0)|}\int_{T_r(x_0)} g(x', t) u(x', 0, t) \varphi(x', 0, t)^2 dx', 
\end{aligned}
\end{equation}
for a.e.~$t \in \Gamma_{\beta, z_0}(r)$, where $dx'$ denotes the $(n-1)$--dimensional Lebesgue measure. Let us consider the last term on the right-hand side of \eqref{energy-1}. By applying the trace inequality and Young's inequality, we infer that
\begin{equation*}\label{bdry-control}
   \begin{aligned}
&\int_{T_r(x_0)}u(x', 0, t)^2 \varphi (x', 0, t)^2dx'\\
&\leq N(n)\left(r\int_{B_r^+(x_0)}|\nabla (u\varphi)(x,t)|^2dx
+r^{-1}\int_{B_r^+(x_0)}(u\varphi)^2(x,t)dx\right).
   \end{aligned}
\end{equation*}
Since $|T_r(x_0)|\leq N(n)r^{-1}|B_r^+(x_0)|$, it follows that, for every $\varepsilon\in(0,1)$,
\begin{equation}\label{bdry-control-2}
\begin{aligned}
&\frac1{|B_r^+(x_0)|}\int_{T_r(x_0)}|g(x',t)u(x', 0,t)|\varphi(x', 0, t)^2dx'\leq\varepsilon\fint_{B_r^+(x_0)}|\nabla(u\varphi)(x,t)|^2dx\\
&\qquad +N(n)r^{-2}\fint_{B_r^+(x_0)}u^2(x,t)dx +\frac{N(n)}{\varepsilon}\fint_{T_r(x_0)}g(x', t)^2dx' \quad \text{for a.e.}\ t \in \Gamma_{\beta, z_0}(r).
\end{aligned}
\end{equation}
From this estimate and by using the ellipticity and boundedness conditions on $\bA$ in \eqref{ellip-cond}, we can follow the standard energy estimates and choose a specific $\varepsilon\in (0,1)$ to obtain
\begin{align*}
& \frac{d}{dt}\fint_{B_r^+(x_0)} \beta(x)u(x,t)^2 \varphi(x,t)^2 dx + \frac{\nu}{4}\fint_{B_r^+(x_0)} |\nabla (u\varphi)(x,t)|^2 dx \\
&\leq N(n,\nu)\left\{\fint_{B_r^+(x_0)} u(x,t)^2\Big(1+r^{-2}+|\nabla \varphi(x,t)|^2 + \beta(x)|\varphi_t(x,t)| \Big) dx \right. \\
&\qquad \left. +\fint_{B_r^+(x_0)} |F(x,t)|^2  dx+\fint_{T_r(x_0)} g(x',t)^2 dx'\right\} \quad \text{for a.e.}\ t \in \Gamma_{\beta, z_0}(r).
\end{align*}
From this, and by integrating the last inequality over $(t_0-r^2\mathcal E_\beta(x_0,r),\tau)$ for $\tau\in\Gamma_{\beta,z_0}(r)$, and taking the essential supremum over $\tau\in \Gamma_{\beta,z_0}(r)$, we obtain the desired estimate.
\end{proof}

Using Lemma~\ref{Caccio-2}, we also have the following result which improves the assertion estimate in Lemma \ref{Caccio-2} by removing some weighted $L^2$ norm of $u(x,t)$ in its right hand side.
\begin{lemma} \label{bdr-improved-Caccio-1} Let $\nu \in (0,1)$ and $M_0 \geq 1$. Suppose that  $u \in \hat{\W}^{1,2}(Q_{2, \beta}^+(z_0), \beta)$ is a weak solution of \eqref{eqn-bdry} with $F \in L^2(Q_{2, \beta}^+(z_0))^n$ and $g \in L^2(T_{2}(x_0) \times \Gamma_{\beta, z_0}(2))$. Then
\begin{align*} \notag
& \sup_{t \in \Gamma_{\beta, z_0}(3/2)}\frac{1}{4\mathcal E_\beta(x_0,2)}\fint_{B_{3/2}^+(x_0)} u(x,t)^2 \beta(x) dx +  \fint_{Q_{3/2, \beta}^+(z_0)} |\nabla u(x,t)|^2dxdt\\ \label{inter-energy-0317-bdry}
& \leq N\fint_{Q_{2,\beta}^+(z_0)} \big( u(x,t)^2 + |F(x,t)|^2\big)dx dt+ N\fint_{T_{2}(x_0)\times\Gamma_{\beta, z_0}(2)} |g(x',t)|^{2} dx'dt,
\end{align*}
where $\mathcal{E}_{\beta}(x_0, r)$ is defined in \eqref{t-height}, and $N = N(n, \nu, M_0)>0$.
\end{lemma}

\begin{proof} As indicated in Remark \ref{A-p-remark}, we have $\bar\beta=\beta^{\frac{n_0}{2}} \in A_{1+\frac{n_0}{2}}$. Then, there is $N=N(n, M_0)>0$ such that
\[
\frac{\bar\beta(B_s(x_0)\setminus B_r(x_0))}{\bar\beta(B_2(x_0))}
\geq N^{-1}\left(\frac{|B_s(x_0)\setminus B_r(x_0)|}{|B_2(x_0)|}\right)^{1+n_0/2}
\geq N^{-1}(s-r)^{1+n_0/2},
\]
for every $3/2\leq r<s\leq2$. Hence, by \eqref{heights},
\[
s^2\mathcal E_\beta(x_0,s)-r^2\mathcal E_\beta(x_0,r)
\geq N^{-1}(s-r)^{1+n_0/2}\mathcal E_\beta(x_0,2).
\]
We can therefore choose a smooth cutoff $0\leq\varphi\leq1$ on $\overline{Q_{2,\beta}(z_0)}$, equal to $1$ in $Q_{r,\beta}(z_0)$ and vanishing outside $Q_{s,\beta}(z_0)$ and near $\partial' Q_{s,\beta}(z_0)$, such that
\[
|\nabla\varphi|\leq\frac{N}{s-r}, \quad \text{and} \quad
|\varphi_t|\leq\frac{N}{(s-r)^{1+n_0/2}\mathcal E_\beta(x_0,2)}.
\]
Applying Lemma~\ref{Caccio-2} on the cylinder $Q_{2,\beta}^+(z_0)$ with this $\varphi$, we obtain
\begin{align*}
&\sup_{t\in\Gamma_{\beta,z_0}(r)}
\frac{1}{4\mathcal E_\beta(x_0,2)}
\fint_{B_r^+(x_0)}\beta(x)u(x,t)^2\,dx
+\fint_{Q_{r,\beta}^+(z_0)}|\nabla u(x,t)|^2\,dxdt\\
&\leq\frac{N}{(s-r)^{1+n_0/2}}\left\{
\fint_{Q_{2,\beta}^+(z_0)}(u(x,t)^2+|F(x,t)|^2)\,dxdt
\right.\\ 
&\hspace{20mm}\left. +\frac{1}{\mathcal E_\beta(x_0,2)}
\fint_{Q_{s,\beta}^+(z_0)}\beta(x)u(x,t) ^2\,dxdt 
+\fint_{T_2(x_0)\times\Gamma_{\beta,z_0}(2)}g(x',t)^2\,dx'dt
\right\},
\end{align*}
where $N=N(n,\nu,M_0)>0$. Here we also used the doubling property to compare the volumes of the cylinders with radii between $3/2$ and $2$. 

\smallskip
Note that $(\beta)_{B_s(x_0)}\leq N\mathcal E_\beta(x_0,2)$ for $3/2\leq s\leq2$. We can apply Lemma~\ref{embedd-lemma} together with Young's inequality to obtain
\begin{align*}
&\frac{1}{(s-r)^{1+n_0/2}\mathcal E_\beta(x_0,2)}
\fint_{Q_{s,\beta}^+(z_0)}\beta(x)u(x,t)^2\,dxdt\\
&\leq\epsilon\fint_{Q_{s,\beta}^+(z_0)}|\nabla u(x,t)|^2\,dxdt
+\frac{N(n,M_0,\epsilon)}{(s-r)^m}
\fint_{Q_{s,\beta}^+(z_0)}u(x,t)^2\,dxdt \qquad \text{for all}\quad \epsilon>0,
\end{align*}
where $m=(1+n_0/2)/(1-\theta)$ and $\theta=\theta(n, M_0)\in(0,1)$ is the constant in Lemma~\ref{embedd-lemma}. Then, combining the last two estimates and choosing $\epsilon$ sufficiently small depending on $n,\nu$ and $M_0$, we infer that
\begin{align*}
&\sup_{t\in\Gamma_{\beta,z_0}(r)}
\frac{1}{4\mathcal E_\beta(x_0,2)}
\fint_{B_r^+(x_0)}\beta(x)u(x,t)^2\,dx
+\fint_{Q_{r,\beta}^+(z_0)}|\nabla u(x,t)|^2\,dxdt\\
&\leq\frac12\fint_{Q_{s,\beta}^+(z_0)}|\nabla u(x,t)|^2\,dxdt+\frac{N}{(s-r)^m}\left\{
\fint_{Q_{2,\beta}^+(z_0)}\left(u(x,t)^2+|F(x,t)|^2\right)dxdt \right.\\
&\qquad \left.
+\fint_{T_2(x_0)\times\Gamma_{\beta,z_0}(2)}g(x',t)^2\,dx'dt
\right\}
\end{align*}
for all $3/2\leq r<s\leq2$, where $N=N(n,\nu,M_0)>0$. We then apply a standard iteration argument, \cite[Lemma 4.3]{Lin} for example, to obtain the assertion in the lemma. The proof of the lemma is complete.
\end{proof}

Our next lemma gives parabolic Poincar\'e type estimates, which are needed to control the local $\hat\W^{1,2}$ norm of weak solutions to \eqref{eqn-bdry}.

\begin{lemma} \label{bdry-Poincare}
Let $\nu\in (0,1)$ and $M_0 \geq 1$. There exists a constant $N = N(n, \nu, M_0)>0$ such that the following assertion holds. Suppose that $r>0$ and $u \in \hat{\W}^{1,2}(Q^+_{r, \beta}, \beta)$ is a weak solution of \eqref{eqn-bdry} in $Q^+_{r, \beta}$. Then
\[
\fint_{Q^+_{r, \beta}} |u(x,t)-(u)_{Q^+_{r, \beta}} |^2 dxdt\leq N r^2\fint_{Q^+_{r, \beta}}\left(|\nabla u(x,t)|^2 + |F(x,t)|^2\right) dxdt.
\]
\end{lemma}

\begin{proof} We follow and adjust the argument in \cite{GM}. By Lemma~\ref{scaling-lemma} and the dilation \eqref{scale-1}, we may assume that $r=1$ and $\mathcal E_\beta(1)=1$. Because of this, we write $Q=Q_1^+=B_1^+\times(-1,0)$, and let $\psi\in C_0^\infty(B_1^+)$ be a nonnegative function such that
\begin{equation}\label{psi-def}
\int_{B_1^+}\psi(x)\,dx=1, \quad \text{and} \quad \|\psi\|_{C^1(B_1^+)}\leq N(n).
\end{equation}
Because $\mathcal E_\beta(1)=1$, and $[\frac{1}{\beta}]_{A_{1+\frac{2}{n_0}}}\leq M_0$,  we apply H\"older's inequality to infer that
\begin{equation}\label{beta prop}
\fint_{B_1}\beta(x)^{-1}\,dx\leq M_0 \quad \text{and} \quad (\beta)_{B_1}\leq1.
\end{equation}

\smallskip
Next, let us denote 
\[
m_0=\int_{B_1^+}\beta(x)\psi(x)\,dx.
\]
It then follows from H\"{o}lder's inequality that
\[
1=\left(\int_{B_1^+}(\beta(x)\psi(x))^{1/2}(\beta(x)^{-1}\psi(x))^{1/2}\,dx\right)^2
\le m_0\int_{B_1^+}\beta(x)^{-1}\psi(x)\,dx.
\]
From this, \eqref{psi-def}, and \eqref{beta prop}, we infer that
\begin{equation}\label{m0-bounds}
N_0^{-1}\leq\Big(\int_{B_1^+}\beta(x)^{-1}\psi(x)\,dx\Big)^{-1}\leq m_0\leq N_0,
\end{equation}
for some $N_0=N_0(n,M_0)\geq 1$.
On the other hand, it follows from Lemma~\ref{embedd-lemma} and Young's inequality that
\begin{equation}\label{L2 beta}
\int_Q\beta(x)u(x,t)^2\,dxdt\leq N(n, M_0)\int_Q\bigl(u(x,t)^2+|\nabla u(x,t)|^2\bigr)\,dxdt<\infty.
\end{equation}

\smallskip
Now, following \cite{GM}, we define 
\[
m(t)=\frac{1}{m_0}\int_{B_1^+}\beta(x)u(x,t)\psi(x)\,dx \qquad \text{for a.e.} \quad t\in(-1,0),
\]
and it is not hard to see that $m\in L^2(-1,0)$ by using H\"older's inequality and \eqref{L2 beta}. On the other hand, by Poincar\'e's inequality, we have
\[
\int_{B_1^+}|u(x,t)-(u)_{B_1^+}(t)|^2\,dx\leq N(n)\int_{B_1^+}|\nabla u(x,t)|^2\,dx, \quad \text{for a.e. } \ t \ \in \ (-1,0).
\]
Thus, applying Lemma~\ref{embedd-lemma} to $u-(u)_{B_1^+}(t)$, and using $(\beta)_{B_1}\leq1$, we obtain
\begin{align*}
\int_Q\beta(x)|u(x,t)-(u)_{B_1^+}(t)|^2\,dxdt\leq N(n, M_0)\int_Q|\nabla u(x,t)|^2\,dxdt.
\end{align*}
Also, by the definition of $m(t)$, H\"older's inequality, \eqref{psi-def} and \eqref{m0-bounds}, we see that
\begin{align*}
|m(t)-(u)_{B_1^+}(t)|^2&\leq\frac{1}{m_0}\int_{B_1^+}\beta(x)\psi(x)|u(x,t)-(u)_{B_1^+}(t)|^2\,dx\\
&\leq N(n, M_0)\int_{B_1^+}\beta(x)|u(x,t)-(u)_{B_1^+}(t)|^2\,dx.
\end{align*}
Then, by combining the last three estimates, we obtain 
\begin{align} \notag
& \int_Q|u(x,t)-m(t)|^2\,dxdt \\ \notag
& \leq2\int_Q|u(x,t)-(u)_{B_1^+}(t)|^2\,dxdt +2|B_1^+|\int_{-1}^0|m(t)-(u)_{B_1^+}(t)|^2\,dt\\ \label{u-m-est}
& \leq N(n, M_0)\int_Q|\nabla u(x,t)|^2\,dxdt.
\end{align}

\smallskip
Next, recalling that $\psi\in C_0^\infty(B_1^+)$, by testing the equation \eqref{eqn-bdry} with $\psi(x)\eta(t)$ with $\eta\in C_0^\infty(-1,0)$, we obtain
\[
-\int_{-1}^0m(t)\eta'(t)\,dt=-\frac{1}{m_0}\int_Q\langle\bA(x,t)\nabla u(x,t)+F(x,t),\nabla\psi(x)\rangle\eta(t)\,dxdt.
\]
This implies that
\[
m'(t)=-\frac{1}{m_0}\int_{B_1^+}\langle\bA(x,t)\nabla u(x,t)+F(x,t),\nabla\psi(x)\rangle\,dx
\]
in the sense of distributions. By \eqref{ellip-cond}, the bound for $\nabla\psi$ in \eqref{psi-def}, and \eqref{m0-bounds}, we have
\[
\int_{-1}^0|m'(t)|^2dt
\leq N(n, \nu, M_0)\int_Q\bigl(|\nabla u(x,t)|^2+|F(x,t)|^2\bigr)\,dxdt.
\]
From this and $m\in L^2(-1,0)$, we see that $m\in W^{1,2}(-1,0)$. Then, applying Poincar\'e's inequality on $(-1,0)$, we obtain
\begin{align}\notag
\int_{-1}^0|m(t)-(m)_{(-1,0)}|^2\,dt
& \leq N\int_{-1}^0|m'(t)|^2\,dt \\ \label{m-m-est}
& \leq N(n, \nu, M_0)\int_Q\bigl(|\nabla u(x,t)|^2+|F(x,t)|^2\bigr)\,dxdt.
\end{align}

\smallskip
Finally, due to the known fact that $(u)_Q$ minimizes $\int_Q|u(x,t)-c|^2\,dxdt$ over all constants~$c$, we take $c=(m)_{(-1,0)}$ and use \eqref{u-m-est} and \eqref{m-m-est} to obtain
\begin{align*}
& \int_Q|u(x,t)-(u)_Q|^2\,dxdt \leq\int_Q|u(x,t)-(m)_{(-1,0)}|^2\,dxdt\\
&\leq2\int_Q|u(x,t)-m(t)|^2\,dxdt
+2|B_1^+|\int_{-1}^0|m(t)-(m)_{(-1,0)}|^2dt\\
&\leq N(n, \nu, M_0)\int_Q\bigl(|\nabla u(x,t)|^2+|F(x,t)|^2\bigr)\,dxdt.
\end{align*}
Dividing this last inequality by $|Q|$, and scaling it back, we obtain the assertion in the lemma. The proof is complete.
\end{proof}

Thanks to Lemma \ref{bdry-Poincare}, we have the following result.
\begin{corollary}\label{bdry-corollary} 
  For every $\nu \in (0,1)$ and $M_0 \geq 1$, there exists a constant $N = N(n, \nu, M_0)>0$ such that the following assertion holds. Suppose that $\beta$ is a weight satisfying \eqref{beta-cond}, $\bA$ satisfies \eqref{ellip-cond} in $Q^+_{1,\beta}$, $F\in L^2(Q^+_{1,\beta})^n$, and $g\in L^2(T_1\times\Gamma_\beta(1))$.
Then, for every weak solution $u \in \hat{\W}^{1,2}(Q^+_{1, \beta}, \beta)$ of \eqref{eqn-bdry} in $Q^+_{1, \beta}$, it holds that
\begin{align*}
\|u-(u)_{Q^+_{1, \beta}} \|_{\hat{\W}^{1,2}(Q^+_{1, \beta}, \beta)} & \leq N \left(\|\nabla u\|_{L^2(Q^+_{1, \beta})} + \|F\|_{L^2(Q^+_{1,\beta})}+\|g\|_{L^2(T_1\times \Gamma_{\beta}(1))}\right).
\end{align*} 
\end{corollary}

\begin{proof} 
By Lemma \ref{bdry-Poincare}, we have
\[
\|u-(u)_{Q_{1,\beta}^+}\|_{L^2(Q_{1,\beta}^+)}
\leq N\left(\|\nabla u\|_{L^2(Q_{1,\beta}^+)}+\|F\|_{L^2(Q_{1,\beta}^+)}\right),
\]
where $N=N(n,\nu,M_0)>0$. On the other hand, for every $\varphi\in \mV_2(B_1^+)$ with $\mV_2(B_1^+)$ defined in \eqref{space-Vp}, the weak formulation gives, for a.e. $t\in\Gamma_\beta(1)$,
\[
\langle\beta u_t(\cdot, t),\varphi\rangle
=-\int_{B_1^+}\langle\bA(x,t)\nabla u(x,t)+F(x,t),\nabla\varphi(x)\rangle\,dx
+\int_{T_1}g(x',t)\varphi(x', 0)\,dx'.
\]
By \eqref{ellip-cond}, the Cauchy--Schwarz inequality, and the trace inequality, it follows that
\begin{align*}
|\langle\beta u_t(t),\varphi\rangle|
\leq N(n,\nu)\left(\|\nabla u(\cdot,t)\|_{L^2(B_1^+)}+\|F(\cdot,t)\|_{L^2(B_1^+)}+\|g(\cdot,t)\|_{L^2(T_1)}\right)\|\varphi\|_{\mV_2(B_1^+)}.
\end{align*}
Due to this, taking the integral over $t\in\Gamma_\beta(1)$, we obtain
\[
\|\beta u_t\|_{L^2(\Gamma_{\beta}(1), \mV_2(B_1^+)^*)}\leq N(n,\nu)\left(\|\nabla u\|_{L^2(Q_{1,\beta}^+)}+\|F\|_{L^2(Q_{1,\beta}^+)}+\|g\|_{L^2(T_1\times \Gamma_{\beta}(1))}\right).
\]
Recalling the definition of the norm in $\hat{\W}^{1,2}(Q^+_{1, \beta}, \beta)$ and combining the estimates above, we obtain the assertion. 
\end{proof}

\subsection{Boundary approximation estimates}\label{bdry approx section} To prove Theorem \ref{inter-theorem-bdry}, we use a perturbation technique by freezing the coefficients. Hence, we consider the following parabolic problem with a homogeneous conormal boundary condition on the flat part $T_{r}(x_0) \times \Gamma_{\beta, z_0}(r)$ as a reference problem for \eqref{eqn-bdry}:
\begin{equation} \label{bdry-approx}
\left\{
\begin{aligned}
(\beta)_{B_r(x_0)} v_t - \textup{div} ((\bA)_{B^+_{r}(x_0)}(t) \nabla v)  &= 0 \quad &&\text{in}\quad Q^+_{r,\beta}(z_0),\\[4pt]
((\bA)_{B^+_{r}(x_0)}(t)\nabla v)\cdot \vec{\mathbf{n}} &=0 \quad &&\text{on} \quad T_{r}(x_0) \times \Gamma_{\beta, z_0}(r).
\end{aligned}\right.
\end{equation}
\smallskip
We first prove the following stability result.
\begin{lemma}\label{bdry-Caccioppoli}
Suppose that $\{\beta_k\}$ is a sequence of weights satisfying \eqref{beta-cond}, $\{\bA_k\}$ is a sequence of matrices satisfying \eqref{ellip-cond} in $Q_4^+$, and $\{F_k\}\subset L^2(Q_4^+)^n$ and $\{g_k\}\subset L^2(T_4\times(-16,0))$. Suppose also that
\begin{equation}\label{conv assumptions}
\begin{aligned}
 &\mathcal E_{\beta_k}(4)=1 \qquad \forall\ k\geq 1,\qquad  \beta_k\to\beta_0 \quad\text{in}\quad L^1(B_4),\\
&\|\bA_k - \bar \bA_k \|_{L^2(Q_4^+)} \to 0,  
\qquad \bar\bA_k\rightharpoonup\bA_0(t)\quad\text{weak--$^*$ in}\quad L^\infty(-16,0),\\
&F_k\to0\quad\text{in}\quad L^2(Q_4^+)^n,
\qquad g_k\to0\quad\text{in}\quad L^2(T_4\times(-16,0))\qquad \text{as}\quad k\to\infty,
\end{aligned}
\end{equation}
for some constant $\beta_0\in [\frac{1}{M_0},1]$ and some matrix $\bA_0(t)$, where $\bar\bA_k=(\bA_k)_{B_4^+}(t)$. Also, let $u_k\in\hat\W^{1,2}(Q_4^+,\beta_k)$ be weak solutions of
\[
\left\{
\begin{aligned}
\beta_k(x) \partial_t u_k -\textup{div}(\bA_k\nabla u_k)&=\textup{div}(F_k)
\quad &&\text{in}\quad Q_4^+,\\
(\bA_k\nabla u_k+F_k)\cdot\vec{\mathbf n}&=g_k
\quad &&\text{on}\quad T_4\times(-16,0),
\end{aligned}\right.
\]
in the sense of Definition~\ref{boundary-weak-def}. If
\begin{equation}\label{conv assumptions-2}
u_k\to u_0\quad\text{in}\quad L^2(Q_4^+),
\qquad \nabla u_k\rightharpoonup\nabla u_0\quad\text{in}\quad L^2(Q_4^+)^n
\end{equation}
for some $u_0$, then $u_0\in\hat\W^{1,2}(Q_4^+)$ is a weak solution of
\begin{equation}\label{limit-eqn}
\left\{
\begin{aligned}
\beta_0 \partial_t u_0 -\textup{div}(\bA_0(t)\nabla u_0)&=0
\quad &&\text{in}\quad Q_4^+,\\
(\bA_0(t)\nabla u_0)\cdot\vec{\mathbf n}&=0
\quad &&\text{on}\quad T_4\times(-16,0).
\end{aligned}\right.
\end{equation}
Moreover, for every $s\in(0,4)$ and $\tau\in(0,16)$, we have
\begin{equation}\label{strong-conv}
\nabla u_k\to\nabla u_0\quad\text{in}\quad L^2(B_s^+\times(-\tau,0))^n.
\end{equation}
\end{lemma}

\begin{proof}
To begin with, we prove some convergence results that will be used in the proof. Note that $\{u_k\}$ is bounded in $L^2((-16,0),W^{1,2}(B_4^+))$, and that $(\beta_k)_{B_4}\leq \mathcal{E}_{\beta_k}(4)=1$ due to H\"older's inequality. Hence, by applying Lemma~\ref{embedd-lemma} to $u_k-u_0$, we obtain
\[
\int_{Q_4^+}\beta_k(x)|u_k(x,t)-u_0(x,t)|^2dxdt\to0.
\]
Also, using $\beta_k\to\beta_0$ in $L^1(B_4)$ and applying Lemma~\ref{embedd-lemma} again, we infer that
\[
\int_{Q_4^+}|\beta_k(x)-\beta_0|u_0(x,t)^2dxdt\to0.
\]
Consequently, from the last two estimates, by using H\"older's inequality and the triangle inequality, we have
\begin{equation}\label{beta-u}
\beta_ku_k\to\beta_0u_0\quad\text{in }L^1(Q_4^+),
\quad \text{and} \quad \beta_ku_k^2\to\beta_0u_0^2\quad\text{in }L^1(Q_4^+).
\end{equation}
Besides, we recall the following trace inequality 
\[
\|v(\cdot, 0)\|_{L^2(T_4)}^2
\leq N \|v\|_{L^2(B_4^+)}\|v\|_{W^{1,2}(B_4^+)},
\qquad v\in W^{1,2}(B_4^+),
\]
where $N = N(n)>0$. Applying this inequality to $v=u_k(\cdot,t)-u_0(\cdot,t)$, and then integrating on $(-16,0)$ and using H\"older's inequality in time, we obtain
\[
\begin{aligned}
& \|u_k(\cdot, 0, \cdot)-u_0(\cdot, 0, \cdot)\|_{L^2(T_4\times(-16,0))}^2 \\
& \leq N \|u_k-u_0\|_{L^2(Q_4^+)}
\|u_k-u_0\|_{L^2((-16,0),W^{1,2}(B_4^+))}\longrightarrow 0, \quad \text{as} \quad k \rightarrow \infty,
\end{aligned}
\]
where the first factor tends to zero due to \eqref{conv assumptions-2} and the second is uniformly bounded due to Corollary \ref{bdry-corollary} and the assumptions \eqref{conv assumptions} and \eqref{conv assumptions-2}. Thus,
\begin{equation}\label{trace-conv}
u_k(\cdot, 0, \cdot) \to u_0(\cdot, 0, \cdot)\quad\text{in}\quad L^2(T_4\times(-16,0)).
\end{equation}
Lastly, as in  \cite[Proof of Lemma~5.7]{FP1}, it can be shown that
\begin{equation}\label{A-conv}
\bA_k\nabla u_k\rightharpoonup\bA_0\nabla u_0,
\qquad \bA_k^T\nabla u_k\rightharpoonup\bA_0^T\nabla u_0
\quad\text{in}\quad L^2(Q_4^+)^n.
\end{equation}

\smallskip
Now, we show that $u_0$ is a solution of \eqref{limit-eqn}. For any fixed $\varphi\in C_0^\infty(B_4\times(-16,0))$, and for every $k\in\mathbb N$, we use $\varphi$ as a test function to the equation of $u_k$  to obtain
\[
\begin{aligned}
& \int_{Q_4^+}\beta_k(x)u_k(x,t)\partial_t\varphi(x,t)\,dxdt
-\int_{Q_4^+}\wei{\bA_k\nabla u_k(x,t),\nabla\varphi(x,t)}dxdt\\
& =\int_{Q_4^+}\wei{F_k(x,t),\nabla\varphi(x,t)}dxdt
-\int_{T_4\times(-16,0)}g_k(x',t)\varphi(x', 0, t)\,dx'dt.
\end{aligned}
\]
We handle term by term in this identity.  For the two terms on its left-hand side, by \eqref{beta-u} and \eqref{A-conv}, it follows that
\[
\lim_{k\to\infty}\int_{Q_4^+}\beta_k(x)u_k(x,t)\partial_t\varphi(x,t)\,dxdt=\int_{Q_4^+}\beta_0u_0(x,t) \partial_t\varphi(x,t)\,dxdt,
\]
and
\[
\lim_{k\to\infty}\int_{Q_4^+}\wei{\bA_k(x,t)\nabla u_k(x,t),\nabla\varphi(x,t)}dxdt
=\int_{Q_4^+}\wei{\bA_0(t)\nabla u_0(x,t),\nabla\varphi(x,t)}dxdt.
\]
On the other hand, the two terms on the right-hand side tend to zero by the convergence of $F_k$ and $g_k$ due to the assumption \ref{conv assumptions}. Combining these, we obtain
\[
\int_{Q_4^+}\beta_0u_0(x,t)\partial_t\varphi(x,t)\,dxdt
-\int_{Q_4^+}\wei{\bA_0(t)\nabla u_0(x,t),\nabla\varphi(x,t)}dxdt=0
\]
for all $\varphi\in C_0^\infty(B_4\times(-16,0))$. This proves $u_0$ is a weak solution of \eqref{limit-eqn}.

\smallskip
Next, we prove \eqref{strong-conv}. Let $s\in(0,4)$ and $\tau\in(0,16)$ be fixed, and let $\chi\in C_0^\infty(B_4)$ and $\eta\in C^\infty([-16,0])$ be cut-off functions such that $0\leq\chi,\eta\leq1$, $\chi=1$ on $B_s$, $\eta=0$ near $-16$, and $\eta=1$ on $[-\tau,0]$. We first claim that
\begin{equation}\label{trace-lower}
\liminf_{k\to\infty}\int_{B_4^+}\beta_k(x)\chi(x)^2u_k(x,0)^2dx\geq
\beta_0\int_{B_4^+}\chi(x)^2u_0(x,0)^2dx.
\end{equation}
We note that due to \cite[Lemma 3.3]{FP1}, the integration terms in \eqref{trace-lower} are well-defined. To prove \eqref{trace-lower}, for every $\psi\in C_0^\infty(B_4)$, let us denote $I=(-16, 0)$ and define
\[
m_k(t)=\int_{B_4^+}\beta_k(x)u_k(x,t)\chi(x)\psi(x)\,dx,\quad \text{for a.e.}~ t\in \overline{I}, \quad k\in \N.
\]
Applying H\"older's inequality to $m_k(t)$ and noting that $(\beta_k)_{B_4}\leq 1$ for all $k\in \N$, we obtain
\begin{equation}\label{L2 of m}
\begin{aligned}
\|m_k\|_{L^2(I)}^2
&\leq \|\chi\psi\|_{L^{\infty}(B_4^+)}^2\left(\int_{B_4^+}\beta_k(x) dx\right)\left(\int_{Q_4^+}\beta_k(x)u_k(x,t)^2 dxdt\right)\\
&\leq \|\chi\psi\|_{L^{\infty}(B_4^+)}^2|B_4|\left(\int_{Q_4^+}\beta_k(x)u_k(x,t)^2 dxdt\right)
\leq N_0 \qquad \text{for all}\quad k\in\N,
\end{aligned}
\end{equation}
for some $N_0>0$ independent of $k$, where in the last inequality we also used \eqref{beta-u}. Besides, testing the equation of $u_k$ with $\chi(x)\psi(x)\rho(t)$ for $\rho\in C_0^{\infty}(I)$, we infer that
\[
m_k'(t)=-\int_{B_4^+}\langle \bA_k\nabla u_k+F_k, \nabla (\chi\psi)\rangle(x,t) dx+\int_{T_4}g_k(x',t)\chi(x',0)\psi(x',0)dx' \quad \text{for a.e.}\quad t\in I
\]
in the sense of distributions. 
Then, by a direct computation, using the fact that $\{u_k\}$ is bounded in $L^2(I, W^{1,2}(B_4^+))$ and the assumptions on $\{F_k\}$ and $\{g_k\}$ in \eqref{conv assumptions}, we have
\begin{align*}
\|m_k'\|_{L^2(I)}
\leq N(n, \nu, \|\chi\psi\|_{C^1(B_4^+)}) \left(\|\nabla u_k\|_{L^2(Q_4^+)}+\|F_k\|_{L^2(Q_4^+)}+\|g_k\|_{L^2(T_4\times I)}\right)\leq N_1
\end{align*}
for some $N_1>0$ independent of $k$. Combining this with \eqref{L2 of m}, we therefore conclude that $\{m_k\}$ is bounded in $H^1(I)$. On the other hand, by \eqref{beta-u}, we have
\[
m_k(t)\rightarrow \beta_0\int_{B_4^+} u_0(x,t)\chi(x)\psi(x)\,dx \quad \text{in} \quad L^1(I), \quad \text{as}\quad k\to\infty.
\] 
Due to this and the compact embedding $H^1(I)\hookrightarrow C(\overline{I})$, it follows that
\begin{equation}\label{limit-m0}
m_k(0)=\int_{B_4^+}\beta_k(x) u_k(x,0)\chi(x)\psi(x)\,dx\longrightarrow
\beta_0\int_{B_4^+} u_0(x,0)\chi(x)\psi(x)\,dx,
\end{equation}
as $k\to\infty$. Moreover, by Young's inequality, we have
\[
\int_{B_4^+}\beta_k(x) u_k(x,0)^2\chi(x)^2dx
\geq
2\int_{B_4^+}\beta_k(x) u_k(x,0)\chi(x)\psi(x)\,dx
-\int_{B_4^+}\beta_k(x)\psi(x)^2dx.
\]
Then, for each fixed $\psi$, using \eqref{limit-m0} and the assumption that $\beta_k\to\beta_0$ in $L^1(B_4)$, we take the lower limit of the last inequality to obtain 
\[
\liminf_{k\rightarrow \infty}\int_{B_4^+}\beta_k(x)u_k(x,0)^2\chi(x)^2dx\geq 2\beta_0\int_{B_4^+} u_0(x,0)\chi(x)\psi(x)\,dx-\int_{B_4^+}\beta_0\psi(x)^2dx.
\]
Due to the density of $C_0^\infty(B_4)$ in $L^2(B_4)$, we can plug $\psi=\chi u_0(x,0)$ in the last estimate to obtain the claim \eqref{trace-lower}.

\smallskip
Now, by using Steklov's average if needed, as in the proof of Lemma~\ref{Caccio-2}, we can test the equation of $u_k$ with $u_k\chi^2\eta^2$ for each $k\geq 1$ to obtain
\begin{align} \notag
&\frac12\int_{B_4^+}\beta_k(x)\chi(x)^2u_k(x,0)^2dx+\int_{Q_4^+}\chi(x)^2\eta(t)^2\langle\bA_k\nabla u_k, \nabla u_k\rangle dxdt\\ \notag
&= \int_{Q_4^+}\beta_k(x)u_k(x,t)^2\chi(x)^2\eta(t)\eta'(t) dxdt
-2\int_{Q_4^+}\chi(x)\eta(t)^2u_k(x,t)\langle\bA_k\nabla u_k,\nabla\chi\rangle dxdt\\ \notag
&\quad-\int_{Q_4^+}\langle F_k,\nabla(u_k\chi^2)\rangle\eta(t)^2dxdt+\int_{T_4\times I}g_k(x', t)u_k(x', 0, t)\chi(x',0)^2\eta(t)^2dx'dt\\ \label{cacci-bdry}
&\longrightarrow \int_{Q_4^+}\beta_0u_0(x,t)^2\chi(x)^2\eta(t)\eta'(t) dxdt
-2\int_{Q_4^+}\chi(x)\eta(t)^2u_0(x,t)\langle\bA_0\nabla u_0,\nabla\chi\rangle dxdt 
\end{align}
as $k\to\infty$.
Here, in \eqref{cacci-bdry}, the convergence of the first term is due to \eqref{beta-u}, the convergence of the second term follows from \eqref{conv assumptions-2} and \eqref{A-conv}, and the last two terms tend to zero due to \eqref{conv assumptions}, the uniform bounds for $u_k$ in $L^2(I, W^{1,2}(B_4^+))$ and the trace inequality. Similarly, as $u_0$ satisfies \eqref{limit-eqn}, we can  test the equation of $u_0$ with $u_0\chi^2\eta^2$ to obtain
\begin{equation}\label{weak u_0}
\begin{aligned}
&\frac12\int_{B_4^+}\beta_0\chi(x)^2u_0(x,0)^2dx+\int_{Q_4^+}\chi(x)^2\eta(t)^2\langle\bA_0\nabla u_0, \nabla u_0\rangle dxdt\\
&= \int_{Q_4^+}\beta_0u_0(x,t)^2\chi(x)^2\eta(t)\eta' (t)dxdt
-2\int_{Q_4^+}\chi(x)\eta(t)^2u_0(x,t)\langle\bA_0\nabla u_0,\nabla\chi\rangle dxdt.
\end{aligned}
\end{equation}
Therefore, comparing \eqref{cacci-bdry} and \eqref{weak u_0}, and using \eqref{trace-lower}, we infer that
\[
\begin{aligned}
\limsup_{k\to\infty}\int_{Q_4^+}\chi(x)^2\eta(t)^2\langle\bA_k\nabla u_k,\nabla u_k\rangle dxdt\leq\int_{Q_4^+}\chi(x)^2\eta(t)^2\langle\bA_0\nabla u_0,\nabla u_0\rangle dxdt.
\end{aligned}
\]
Furthermore, since $\{\bA_k\}$ satisfies \eqref{ellip-cond} in $Q_4^+$,  it follows from \eqref{conv assumptions} and \eqref{A-conv} that
\[
\begin{aligned}
& \limsup_{k\to\infty}\int_{Q_4^+}\chi(x)^2\eta(t)^2|\nabla u_k(x,t)-\nabla u_0(x,t)|^2 dxdt\\
& \leq\frac{1}{\nu}\limsup_{k\to \infty}\int_{Q_4^+}\chi(x)^2\eta(t)^2\langle\bA_k(\nabla u_k-\nabla u_0),\nabla u_k-\nabla u_0\rangle dxdt\\
&=\frac{1}{\nu}\limsup_{k\to \infty}\int_{Q_4^+}\chi(x)^2\eta(t)^2\Big[\langle\bA_k\nabla u_k,\nabla u_k\rangle
-\langle\bA_0\nabla u_0,\nabla u_0\rangle\Big] dxdt.
\end{aligned}
\]
Hence, from the last two estimates, we conclude that
\[
\limsup_{k\to\infty}\int_{Q_4^+}\chi(x)^2\eta(t)^2|\nabla u_k(x,t)-\nabla u_0(x,t)|^2dxdt\leq0.
\]
Noting that $\chi=1$ on $B_s$ and $\eta=1$ on $[-\tau,0]$, we obtain \eqref{strong-conv}. Therefore, the lemma is proved.
\end{proof}

The following proposition provides the closeness of solutions of \eqref{eqn-bdry} and \eqref{bdry-approx}.
\begin{proposition}\label{L2-comparision-bdry}
Let $\nu\in(0,1)$ and $M_0\geq1$. For every $\epsilon\in(0,1)$, there exists  a constant $\bar\kappa=\bar\kappa(n,\nu,M_0,\epsilon)\in (0,1)$ such that the following assertion holds. Suppose that $\beta$ satisfies \eqref{beta-cond}, $\bA$ satisfies \eqref{ellip-cond} in $Q_{4,\beta}^+$, and
\begin{equation}\label{smallness hypothesis}
   \begin{aligned}
& \fint_{Q_{4,\beta}^+}\Theta_{\bA,4}(x,t)^2dxdt +\frac1{\beta(B_4)}\int_{B_4}\Theta_{\beta,4}(x)^2dx\\
& \qquad \qquad +\fint_{Q_{4,\beta}^+}|F(x,t)|^2dxdt
+\fint_{T_4\times\Gamma_\beta(4)}g(x',t)^2dx'dt\leq\bar\kappa^2.
\end{aligned}
\end{equation}
Then, for every weak solution $u\in\hat\W^{1,2}(Q_{4,\beta}^+,\beta)$ of \eqref{eqn-bdry} with $\fint_{Q_{4,\beta}^+}|\nabla u(x,t)|^2dxdt\leq1$, there is a weak solution $v\in\hat\W^{1,2}(Q_{4,\beta}^+)$ of \eqref{bdry-approx} such that
\begin{equation}\label{L-2-w-small-bdry}
\fint_{Q_{4,\beta}^+}|u(x,t)-v(x,t)|^2dxdt\leq\epsilon^2
\quad \text{and} \quad \fint_{Q_{2,\beta}^+}|\nabla u(x,t)-\nabla v(x,t)|^2dxdt\leq\epsilon^2.
\end{equation}
Moreover, there exists $N=N(n,\nu,M_0)>0$ such that 
\begin{equation}\label{L-2-mu-bdry}
\fint_{Q_{2,\beta}^+}|\nabla v(x,t)|^2dxdt\leq N,
\quad \text{and} \quad \frac{1}{(\beta)_{B_2}}\fint_{Q_{2,\beta}^+}\beta(x)|u(x,t)-v(x,t)|^2dxdt
\leq N\epsilon^2.
\end{equation}
\end{proposition}

\begin{proof}
By applying the time dilation \eqref{scale-2} and subtracting a constant from $u$, we may assume $\mathcal E_\beta(4)=1$ and $(u)_{Q_4^+}=0$. As such, we can simply write 
\[
Q_{4,\beta}^+=Q_4^+ \quad \text{and}\quad  \Gamma_\beta(4)=(-16,0)
\] 
in the rest of the proof. We start to prove \eqref{L-2-w-small-bdry} by a contradiction argument. Suppose that \eqref{L-2-w-small-bdry} fails for some $\epsilon_0\in (0, 1)$. That is, there are sequences $\{\beta_k\},\{\bA_k\}, \{F_k\},\{g_k\},\{u_k\}$ satisfying the hypotheses in this proposition with $\bar\kappa=1/k\to0$, while for all $k\geq 1$, every weak solution $v$ to the problem 
\begin{equation}\label{bdry-approx-2}
\left\{
\begin{aligned}
(\beta_k)_{B_4} \partial_t v - \textup{div} ((\bA_k)_{B^+_4}(t) \nabla v)  &= 0 \quad &&\text{in}\quad Q^+_{4,\beta_k},\\[4pt]
((\bA_k)_{B^+_{4}}(t)\nabla v)\cdot \vec{\mathbf{n}} &=0 \quad &&\text{on} \quad T_{4} \times \Gamma_{\beta_k}(4)
\end{aligned}\right.
\end{equation}
satisfies either
\begin{equation}\label{L2-contradiction}
\fint_{Q_{4}^+}|u_k(x,t)-v(x,t)|^2dxdt>\epsilon_0^2 \
\quad \ \text{or}\quad \fint_{Q_{2,\beta_k}^+}|\nabla u_k(x,t)-\nabla v(x,t)|^2dxdt>\epsilon_0^2.
\end{equation}
Here again, note that without loss of generality, we assume that $\mathcal{E}_{\beta_k}(4)=1$ and $(u_k)_{Q_4^+}=0$ for all $k\geq 1$. It follows from Corollary~\ref{bdry-corollary} and the assumptions in the proposition that
\[
\|u_k\|_{\hat \W^{1,2}(Q_{4}^+, \beta_k)}\leq N(n, \nu, M_0), \qquad \text{for all}\quad k\geq 1.
\] 
Then, by Theorem \ref{compactness-lemma}, there exists $u_0\in\hat\W^{1,2}(Q_4^+)$ such that, up to a subsequence,
\begin{equation}\label{u-conv}
u_k\to u_0 \quad \text{in}\quad L^2(Q_4^+),\qquad
\nabla u_k\rightharpoonup\nabla u_0 \quad\text{in}\quad L^2(Q_4^+) \qquad \text{as}\quad k\rightarrow \infty.
\end{equation}
Next, applying Lemma \ref{beta-seq} to $\{\beta_k\}$ and noting that $\{\bA_k\}$ is uniformly bounded, we have that, up to subsequences,
\begin{equation}\label{beta-A-conv}
\left\{ 
\begin{split}
& (\beta_k)_{B_4}\to\beta_0 \quad \text{in} \ \mathbb{R}, \quad
\beta_k\to\beta_0 \text{ in }L^1(B_4),\\
& (\bA_k)_{B_4^+}(t)\rightharpoonup\bA_0(t) \quad \text{ weak--$^*$ in } L^\infty((-16, 0)), \quad \text{as} \quad k\rightarrow \infty,
\end{split} \right. 
\end{equation}
for some constant $\beta_0 \in [\frac{1}{M_0},1]$ and some matrix $\bA_0(t)$ satisfying \eqref{ellip-cond}. Also, due to the hypotheses with $\bar\kappa=1/k$ and $k\rightarrow \infty$,  we have 
\[
\|\bA_k-(\bA_k)_{B_4^+}\|_{L^2(Q_4^+)}\to 0,\quad F_k\to0 \quad \text{in}\quad L^2(Q_4^+),\quad g_k\to0 \quad \text{in}\quad L^2(T_4\times(-16,0)).
\]
Therefore, from the last three formulas, we can apply Lemma~\ref{bdry-Caccioppoli} to conclude that $u_0$ solves
\begin{equation} \label{u-0-eqn}
   \left\{
   \begin{aligned}
\beta_0 \partial_t u_0  - \textup{div} (\bA_0(t) \nabla u_0) &= 0 \quad &&\text{in} \quad Q_4^+,\\
(\bA_0(t) \nabla u_0)\cdot \vec{\mathbf{n}} &=0 \quad &&\text{on} \quad T_4 \times (-16, 0).
   \end{aligned}\right.
\end{equation}
Moreover, as $k\rightarrow \infty$, we have that
\begin{equation}\label{grad-u-conv}
\nabla u_k\to\nabla u_0\quad  \text{in}\quad L^2(B_s^+\times(-\tau,0)) \quad \text{for every}\ 0<s<4 \ \text{ and }\ 0<\tau<16.
\end{equation}

\smallskip
Now, we construct a sequence of solutions $\{v_k\}$ of problem \eqref{bdry-approx-2}. Let
$\bar\bA_k(t)=(\bA_k)_{B_4^+}(t)$. Since $u_0\in \hat\W^{1,2}(Q_4^+)$ solves \eqref{u-0-eqn}, we have that $(u_0)_t\in L^2\bigl((-16,0),\mV_2(B_4^+)^*\bigr)$. Note also that $M_0^{-1}\leq (\beta_k)_{B_4}\leq1$ for all $k\in \N$ due to $\mathcal{E}_{\beta_k}(4)=1$. Then, it follows from the classical existence theory (see \cite[Theorem~1.1]{Fackler}, for example) that there exists a unique weak solution 
\[
w_k\in L^2\bigl((-16,0), \mV_2(B_4^+)\bigr)\cap H^1\bigl((-16,0),\mV_2(B_4^+)^*\bigr)
\]
solving the equation
\begin{equation}\label{w-k-eqn}
\left\{
\begin{aligned}
(\beta_k)_{B_4}\partial_t w_k -\textup{div}(\bar{\mathbf A}_k(t)\nabla w_k)
&=-(\beta_k)_{B_4}\partial_t u_0 +\textup{div}(\bar{\mathbf A}_k(t)\nabla u_0) \quad \text{in}\quad Q_4^+,\\
\bigl(\bar{\mathbf A}_k(t)\nabla(w_k+u_0)\bigr)\cdot\vec{\mathbf n}&=0 \qquad \qquad \text{on}\quad T_4\times(-16,0),\\
w_k&=0 \qquad \qquad \text{on}\quad (\partial B_4\cap\mathbb R^n_+)\times(-16,0),\\
w_k(\cdot,-16)&=0\qquad \qquad \text{in}\quad B_4^+.
\end{aligned}
\right.
\end{equation}
Let $v_k=u_0+w_k$ for each $k\in \N$. Then $v_k$ is a weak solution of the corresponding frozen equation \eqref{bdry-approx-2}. Moreover, the standard energy estimate implies that
\[
\|v_k\|_{L^2((-16,0),W^{1,2}(B_4^+))}
+\|(v_k)_t\|_{L^2((-16,0),\mV_2(B_4^+)^*)}
\leq N(n,\nu,M_0).
\]

\smallskip
Next, we apply the classical Aubin--Lions theorem (see \cite[Theorem~1]{Simon}, for example) to $\{v_k\}$  to obtain a function $\bar v\in L^2((-16,0), W^{1,2}(B_4^+))$ with $\bar v_t\in L^2((-16,0), \mV_2(B_4^+)^*)$, such that, up to a subsequence,
\begin{equation}\label{v-k-conv}
\left\{
\begin{aligned}
v_k\to \bar v\quad&\text{in}\quad L^2(Q_4^+),\\
\nabla v_k\rightharpoonup\nabla \bar v\quad&\text{in}\quad L^2(Q_4^+),\\
\partial_t v_k \rightharpoonup \bar v_t\quad&\text{in}\quad L^2((-16,0), \mV_2(B_4^+)^*).
\end{aligned}\right.
\end{equation}
From this, \eqref{beta-A-conv}, and the bounds $M_0^{-1}\leq (\beta_k)_{B_4}\leq 1$, we can follow the proof of Lemma~\ref{bdry-Caccioppoli} to show that $\bar v$ satisfies the same equation as $u_0$ in \eqref{u-0-eqn}, and that
\begin{equation}\label{grad-v-conv}
\nabla v_k\rightarrow \nabla \bar v \quad \text{in}\quad L^2(B_s^+\times(-\tau,0)) \quad \text{for every}\ 0<s<4 \ \text{ and }\ 0<\tau<16.
\end{equation}
Moreover, for $w_k=v_k-u_0$, \eqref{v-k-conv} also implies that
\[
w_k\rightharpoonup  \bar v-u_0\quad\text{in}\quad L^2((-16, 0), \mV_2(B_4^+))\cap H^1((-16,0), \mV_2(B_4^+)^*)
\]
as $k\rightarrow \infty$.  In particular, it follows from \eqref{w-k-eqn} that 
\[
\bar v=u_0 \quad  \text{on}\quad  (\partial B_4\cap\mathbb R^n_+)\times (-16,0)
\] 
in the sense of traces. Due to the continuous embedding from $L^2((-16, 0), \mV_2(B_4^+))\cap H^1((-16,0), \mV_2(B_4^+)^*)$ to $C([-16,0], L^2(B_4^+))$ and the condition that $w_k(\cdot, -16)=0$ in $B_4^+$ for all $k\in \N$ in \eqref{w-k-eqn}, we also have that 
\[
(\bar v-u_0)(\cdot,-16)=0 \quad \text{in} \quad L^2(B_4^+).
\]
From the last two observations and the fact that $u_0$ and $\bar v$ satisfy the same equation  \eqref{u-0-eqn}, by uniqueness, we conclude that $\bar v=u_0$. Consequently, using \eqref{u-conv} and \eqref{v-k-conv}, we have
\begin{equation}\label{conv-1}
\|u_k-v_k\|_{L^2(Q_4^+)}\leq \|u_k-u_0\|_{L^2(Q_4^+)}+\|v_k-\bar v\|_{L^2(Q_4^+)}\rightarrow 0 \quad \text{as}\quad k\rightarrow\infty.
\end{equation}
On the other hand, since $\{\beta_k\}$ satisfies \eqref{beta-cond} and $\mathcal E_{\beta_k}(4)=1$ for each $k$, by the definition of weighted cylinders in \eqref{half cylinder} and the properties of $\beta_k$ in Proposition \ref{Ap-prop}, it follows that
\[
|Q_{2, \beta_k}^+|\geq N_0 \quad \text{for all}\quad k\geq 1, \quad \text{for some} \quad N_0=N_0(n, M_0)>0.
\]
Moreover, we can take $Q^+=B_s^+\times(-\tau,0)$ with $s\in (15/4,4)$ and $\tau\in (15, 16)$ so that
\[
Q_{2, \beta_k}^+\subset Q^+ \quad  \text{for all}\quad k\geq 1.
\]
Thus, using \eqref{grad-u-conv}, \eqref{grad-v-conv} and the fact that $\bar v=u_0$, we can also conclude that
\begin{equation}\label{conv-2}
   \begin{aligned}
&\fint_{Q_{2,\beta_k}^+}|\nabla u_k(x,t)-\nabla v_k(x,t)|^2\,dxdt\\
&\leq\frac{2}{N_0}\|\nabla u_k-\nabla u_0\|_{L^2(Q^+)}^2+\frac{2}{N_0}\|\nabla v_k-\nabla \bar v\|_{L^2(Q^+)}^2\rightarrow 0 \quad \text{as}\quad k\rightarrow\infty.
   \end{aligned}
\end{equation}
However, both convergences \eqref{conv-1} and \eqref{conv-2} contradict \eqref{L2-contradiction}. Therefore, \eqref{L-2-w-small-bdry} is proved.

\smallskip
Lastly, we prove \eqref{L-2-mu-bdry}. Since we assumed that $(u)_{Q_4^+}=0$ at the beginning of the proof, it follows from Corollary \ref{bdry-corollary} and the assumptions in the proposition that
\begin{align*}
\fint_{Q^+_{4}} u(x,t)^2 dx dt  
&\leq N \fint_{Q^+_{4}} \left(|\nabla u(x,t)|^2 + |F(x,t)|^2\right) dxdt+N\fint_{T_4\times\Gamma_{\beta}(4)}g(x', t)^2dx'dt\\
&\leq N\big(1 + \bar\kappa^2) \leq 2N,
\end{align*}
where $N = N(n, \nu, M_0)>0$.  By the last estimate, \eqref{L-2-w-small-bdry}, and the Cauchy-Schwarz inequality, we see that
\begin{align*}
\fint_{Q^+_{4}}v(x,t)^2 dx dt 
\leq 2\fint_{Q^+_{4}} \left(u(x,t)^2+[u(x,t)-v(x,t)]^2\right)dxdt\leq 4N+2\epsilon^2 \leq  4N+2.
\end{align*}
From this, and by the standard energy estimate for the equation of $v$ and the doubling properties of $\beta$ stated in Proposition~\ref{Ap-prop}, we obtain
\[
\fint_{Q^+_{7/2, \beta}} |\nabla v(x,t)|^2 dx dt \leq N \fint_{Q^+_{4}} v(x,t)^2 dx dt \leq N(n, \nu, M_0).
\]
This proves the first assertion in \eqref{L-2-mu-bdry}. For the second one, due to the first conclusion of \eqref{L-2-w-small-bdry} and the doubling properties of $\beta$, it follows that
\[
\fint_{Q_{2, \beta}^+}|u(x,t)-v(x,t)|^2dxdt\leq N(n, M_0)\fint_{Q_4^+}|u(x,t)-v(x,t)|^2dxdt\leq N(n, M_0)\epsilon^2.
\]
Using this and the second conclusion of \eqref{L-2-w-small-bdry}, we can apply Lemma \ref{embedd-lemma} to obtain that
\[
\begin{aligned}
\frac{1}{(\beta)_{B_{2}}}\fint_{Q^+_{2,\beta}} &\beta(x)|u(x,t)-v(x,t)|^2 dx dt
\leq N\Bigl(\fint_{Q^+_{2,\beta}} |u(x,t)-v(x,t)|^2 dxdt\Bigr)^{1-\theta}\\
&\times \Bigl[ \Bigl(\fint_{Q^+_{2, \beta}} |u(x,t)-v(x,t)|^{2} dxdt \Bigr)^{\theta} +  \Bigl(\fint_{Q^+_{2,\beta}} |\nabla (u(x,t)-v(x,t))|^{2} dxdt \Bigr)^{\theta} \Bigr]\\
&\leq N(n, \nu, M_0) \epsilon^2.
\end{aligned}
\]
Therefore, the second assertion of \eqref{L-2-mu-bdry} is also proved.
\end{proof}

\smallskip
Once Proposition \ref{L2-comparision-bdry} is established, we can follow a similar argument as in \cite[Section~6.3]{FP1} and \cite{Phan-2} to derive level set estimates, and then prove Theorem~\ref{inter-theorem-bdry} by adapting the proof of \cite[Theorem~6.1]{FP1}. Since these are standard, we omit the details here.


\section{Global Regularity Estimates and Proof of Theorem \ref{main-theorem}} \label{global-section}
In this section, we prove Theorems \ref{main-thm-2} and \ref{main-theorem}.  We refer the reader to Definitions~\ref{W-space} and \ref{interior-weak-def} for the relevant solution spaces and weak formulations, respectively.  The proof of Theorem~\ref{main-thm-2} closely follows that of \cite[Theorem~7.1]{FP1}, with the necessary modifications to handle the boundary condition. For completeness, we provide the proof here.

\begin{proof}[Proof of Theorem \ref{main-thm-2}]
For the given $\nu \in (0,1)$, $M_0 \geq 1$, and $p \in [2, \infty)$, let 
\begin{equation}\label{choice of delta}
\bar\delta=\min\left\{\delta_0 (n, \nu, p, M_0),\, \frac{\delta_1(n, \nu, p, M_0)}{N_1} \right\}\in (0,1),
\end{equation}
where $\delta_0$ is the constant from the interior regularity estimate in \cite[Theorem~5.1]{FP1},  $\delta_1$ is the constant defined in Theorem~\ref{inter-theorem-bdry},   and $N_1=N_1(n, \nu, M_0)>0$ is the constant given in \eqref{matrix flattening} below. We prove the theorem under this choice of $\bar\delta$ in \eqref{choice of delta}.

\smallskip
As $\partial\Omega\in C^1$, there exists $R=R(\bar\delta,\Omega,R_0)\in(0,R_0)$ sufficiently small
such that the following assertion holds: for each $x_0=(x'_0, x_0^n)\in \partial \Omega$ with $ x_0' \in \R^{n-1}$, there exists a Lipschitz continuous function $\phi: \R^{n-1} \rightarrow \R$ such that, upon relabeling and reorienting the coordinate axes,
\[
\Omega\cap B_R(x_0)=\{(x',x^n)\in B_R(x_0):x^n>\phi(x')\},\qquad\phi(x'_0)=x_0^n,\qquad\nabla\phi(x'_0)=0,
\]
and 
\begin{equation*}\label{h-lip}
\|\nabla\phi\|_{L^\infty(\R^{n-1})}\leq \bar\delta.
\end{equation*}
Taking a linear translation, we can assume without loss of generality that $x_0=0$. Let us define the flattening map $\Phi: \R^{n}\rightarrow \R^n$ as 
\[
\Phi(x',x^n)=(x', x^n-\phi(x')), \quad\text{where}\quad (x', x^n) \in \R^{n-1} \times \R;
\]
and its inverse map $\Phi^{-1}$ is also well-defined, given by
\[
\Phi^{-1}(y',y^n)=(y', y^n+\phi(y')), \quad\text{where}\quad  (y', y^n) \in \R^{n-1} \times \R.
\] 
Note that
\begin{equation*}
\nabla \Phi(x)= \begin{bmatrix}
I_{n-1}& 0 \\
-(\nabla \phi)^* &1
\end{bmatrix}(x),
\quad \text{and}\quad 
\nabla \Phi^{-1}(y)= \begin{bmatrix}
I_{n-1}& 0 \\
(\nabla \phi)^* &1
\end{bmatrix}(y),
\end{equation*}
where $(\nabla \phi)^*$ denotes the transpose matrix of the matrix $\nabla \phi$. As a consequence, due to $\bar\delta\in (0,1)$, 
for each $y=(y', y^n)\in B_{r}(y_0)$ with any $y_0\in \R^n$ and any $r>0$, we have
\begin{align*}
|\Phi^{-1}(y)-\Phi^{-1}(y_0)|
&=\sqrt{|y'-y'_0|^2+|y^n+\phi(y')-y^n_0-\phi(y'_0)|^2}\\
&\leq \sqrt{|y'-y'_0|^2+(|y^n-y^n_0|+\bar\delta r)^2}\\
&\leq \sqrt{r^2+2\bar\delta r^2+\bar\delta^2r^2}=(1+\bar\delta)r<2r.
\end{align*}
Similarly, for each $x=(x', x^n)\in B_{r/2}(\Phi^{-1}(y_0))$, we have
\begin{align*}
|\Phi(x)-y_0|
&=|(x', x^n-\phi(x'))-(y'_0, y^n_0)|\\
&\leq \sqrt{|x'-y'_0|^2+|x^n-\phi(x')-(y^n_0+\phi(y'_0))+\phi(y'_0)|^2}\\
&\leq \sqrt{|x'-y'_0|^2+\Big(|x^n-(y^n_0+\phi(y'_0))|+\frac{\bar\delta r}{2}\Big)^2}\leq r.
\end{align*}
Besides, similarly to \cite[Proof of Theorem~7.1]{FP1}, by a direct computation, we have the following inclusion property: for any $r>0$ and $y_0\in \R^n$, 
\begin{equation}\label{phi-lip-1}
 B_{r/2}(\Phi^{-1}(y_0))\subset  \Phi^{-1}(B_{r}(y_0))\subset B_{2r}(\Phi^{-1}(y_0)).
\end{equation}

\smallskip
Next, we set
\begin{equation*}
\begin{aligned}
&\tilde u(y,t)=u(\Phi^{-1}(y), t),\quad  \tilde \beta(y)=\beta(\Phi^{-1}(y)),\\
&\tilde \bA(y,t)=\nabla\Phi(\Phi^{-1}(y))\bA(\Phi^{-1}(y),t)\nabla \Phi(\Phi^{-1}(y))^*,\\ 
& \tilde F(y,t)=\nabla\Phi(\Phi^{-1}(y))F(\Phi^{-1}(y),t), \quad \tilde g((y',0),t)=g((y',\phi(y')),t)
\sqrt{1+|\nabla\phi(y')|^2},
\end{aligned}
\end{equation*}
where $\nabla \Phi(\Phi^{-1}(y))^*$ denotes the transpose of the matrix $\nabla \Phi(\Phi^{-1}(y))$.  Moreover, we fix an arbitrary $\tilde x\in \Omega$ and denote $B=B_{\bar R}(\tilde x)$ with $\bar R=\textup{diam}(\Omega)+R_0$. It then follows from the doubling property of $\beta^{\frac{n_0}{2}}$ and the definition of $\lambda_0$ in \eqref{lambda-zero-def} that
\[
\left(\int_B\beta(x)^{\frac{n_0}{2}}dx\right)^{2/n_0}\leq N\left(\int_{\Omega}\beta(x)^{\frac{n_0}{2}} \,dx\right)^{2/n_0}=N\lambda_0,
\]
where $N=N(n,M_0,\Omega,R_0)>0$. Applying \textup{(v)} of Proposition~\ref{Ap-prop} to $\beta^{\frac{n_0}{2}}$ on $B$ and using \eqref{phi-lip-1}, we get
\[
r^2\mathcal E_{\tilde\beta}(0,r)\leq N\lambda_0 r^{\kappa}\quad \text{and}\quad
r^2\mathcal E_{\beta}(x,r)
\leq N \lambda_0 r^{\kappa}\qquad \text{for all}\quad x\in\overline\Omega,\quad r\in (0, R/2),
\]
where $\kappa=\kappa(n, M_0)>0$. Thus, we can choose a $\rho=\rho(n,M_0,\Omega,R_0, R, T_0/\lambda_0)>0$ small enough so that
\[
4\Lambda\rho<R,\quad (2\Lambda\rho)^2\mathcal E_{\tilde\beta}(0, 2\Lambda\rho)\leq T_0/2,\quad \text{and}\quad (2\Lambda\rho)^2\mathcal E_{\beta}(x, 2\Lambda\rho)\leq T_0/2 \quad \text{for all}\quad x\in\overline\Omega.
\]
As such, for any $z_0=(0,t_0)$ with $t_0\in[T_0,T]$, we have
\[
B_{2\Lambda\rho}^+\subset\Phi(\Omega\cap B_R),\qquad
Q_{2\Lambda\rho,\tilde\beta}^+(z_0)
\subset B_{2\Lambda\rho}^+\times(0,t_0],
\]
and
\begin{equation}\label{map subset}
\Phi^{-1}(B_{2\Lambda\rho}^+)\times
\Gamma_{\tilde\beta,z_0}(2\Lambda\rho)
\subset(\Omega\cap B_R)\times(0,t_0].
\end{equation}
Moreover, for $x=(x',\phi(x'))\in\partial\Omega\cap B_R$, the outward unit normal $\vec{\mathbf n}(x)$ satisfies
\[
\vec{\mathbf n}(x)=\frac{(\nabla\phi(x'),-1)^*}{\sqrt{1+|\nabla\phi(x')|^2}},
\qquad
[\nabla\Phi(x)]^*(-\vec{\boldsymbol e}_n)
=\sqrt{1+|\nabla\phi(x')|^2}\vec{\mathbf n}(x).
\]
Note that $dS_x=\sqrt{1+|\nabla\phi(x')|^2}\,dx'$. By changing variables in the weak
formulation, it follows that $\tilde u$ solves
\[
\left\{\begin{aligned}
\tilde\beta\tilde u_t-\textup{div}(\tilde\bA\nabla\tilde u)
&=\textup{div}(\tilde F)
&&\text{in }Q_{2\Lambda\rho,\tilde\beta}^+(z_0),\\
(\tilde\bA\nabla\tilde u+\tilde F)\cdot(-\vec{\boldsymbol e}_n)
&=\tilde g
&&\text{on }T_{2\Lambda\rho}\times
\Gamma_{\tilde\beta,z_0}(2\Lambda\rho).
\end{aligned}\right.
\]
On the other hand, similar to \cite[Lemma~7.2]{FP1}, we have that $\tilde\bA$ satisfies \eqref{ellip-cond}, and that $\tilde \beta$ satisfies \eqref{beta-cond} and
\begin{equation}\label{matrix flattening}
\frac{1}{\tilde\beta(B_r(y_0))}
\int_{B_r(y_0)}
|\tilde\beta(y)-(\tilde\beta)_{B_r(y_0)}|^2
\tilde\beta(y)^{-1}\,dy\leq N_1^2\bar\delta^2
\end{equation}
for all $y_0\in B_\rho^+$ and $0<r<\rho$, where $N_1=N_1(n, \nu, M_0)>0$.

\smallskip
Recall the choice of $\bar\delta$ in \eqref{choice of delta}.  After a translation and dilation, we can apply the gradient estimates in Theorem \ref{inter-theorem-bdry} to get
\begin{align*}
\|\nabla\tilde u\|_{L^p(Q_{\rho,\tilde\beta}^+(z_0))}
\leq N\Big[&\rho^{-1}|Q_{\rho,\tilde\beta}^+(z_0)|^{1/p-1/2}\|\tilde u\|_{L^2(Q_{2\Lambda\rho,\tilde\beta}^+(z_0))}\\
&+\|\tilde F\|_{L^p(Q_{2\Lambda\rho,\tilde\beta}^+(z_0))}+\rho^{1/p}\|\tilde g\|_{L^p(T_{2\Lambda\rho}\times\Gamma_{\tilde\beta,z_0}(2\Lambda\rho))}\Big],
\end{align*}
where $N=N(n,\nu,p,M_0)$. By \eqref{phi-lip-1}, we have
\[
\Omega\cap B_{\rho/2}\subset\Phi^{-1}(B_\rho^+),\qquad
\Omega_T\cap Q_{\rho/2,\beta}(z_0)\subset\Phi^{-1}(B_\rho^+)\times\Gamma_{\tilde\beta,z_0}(\rho).
\]
Recall the definition of $\lambda_0$ in \eqref{lambda-zero-def}. Since $p\geq2$, the doubling property of $\beta^{n_0/2}$ and a change of variables imply that
\[
\rho^{-1}|Q_{\rho,\tilde\beta}^+(z_0)|^{1/p-1/2} \leq \bar N(n,M_0,\Omega,\rho)\lambda_0^{1/p-1/2}.
\]
Note also that $|\nabla u|\leq2|\nabla\tilde u|$ and $|\det(\nabla \Phi)|=|\det(\nabla \Phi^{-1})|=1$. Thus, from the last three estimates, we obtain
\begin{equation}\label{bdry esti}
\begin{aligned}
\|\nabla u\|_{L^p(\Omega_T\cap Q_{\rho/2,\beta}(z_0))}
\leq N\Big[&
\lambda_0^{1/p-1/2}
\|u\|_{L^2((\Omega\cap B_R)\times(0,T))}\\
&+\|F\|_{L^p((\Omega\cap B_R)\times(0,T))}
+\rho^{1/p}\|g\|_{L^p((\partial\Omega\cap B_R)\times(0,T))}
\Big],
\end{aligned}
\end{equation}
where $N=N(n,\nu,p,M_0,\Omega,\rho)>0$.

\smallskip
Next, we cover a layer of the parabolic boundary $\partial'(\overline\Omega\times[T_0,T])$ by cylinders with
$\rho_i=\rho/2$, and the remaining interior region by cylinders
with $\rho_i=\rho/(32\Lambda)$. By compactness of $\overline\Omega\times[T_0,T]$ and the doubling
property of $\beta^{n_0/2}$, we obtain a finite covering
\[
\overline{\Omega}\times[T_0,T]
\subset \bigcup_{i=1}^m Q_{\rho_i/2,\beta}(z_i),
\]
where $m\leq N(n,\nu,p,M_0,R_0,\Omega,T_0/\lambda_0,T/\lambda_0)$, and the cylinders can be chosen so that, for each $i$, either
\[
Q_{2\Lambda\rho_i,\beta}(z_i)\subset\Omega\times(0,T]
\]
or
\[
z_i=(x_i,t_i)\in\partial\Omega\times(T_0,T]
\quad\text{with}\quad
\Gamma_{\beta,z_i}(2\Lambda\rho_i)\subset(0,T].
\]
For the former case, we can apply the interior regularity estimates, \cite[Theorem~5.1]{FP1}, to $u$. For the latter case, after taking a linear translation and a dilation, we can apply \eqref{bdry esti} on it. Thus,
\[
\|\nabla u\|_{L^p(\Omega\times(T_0,T))}\leq N\left(\lambda_0^{1/p-1/2}\|u\|_{L^2(\Omega_T)}
+\|F\|_{L^p(\Omega_T)}+\|g\|_{L^p(\partial\Omega\times(0,T))}\right),
\]
where $N=N(n,\nu,p,M_0,R_0,\Omega,T_0/\lambda_0,T/\lambda_0)>0$. By the equation \eqref{main-eqn-2} and the trace theorem, we can also obtain 
\[
\|\beta u_t\|_{L^p((T_0,T), W^{1,p'}(\Omega)^{*})}\leq N\left(\|\nabla u\|_{L^p(\Omega\times(T_0,T))}+\|F\|_{L^p(\Omega_T)}+\|g\|_{L^p(\partial\Omega\times(0,T))}\right).
\]
Finally, the same estimate for $\|u\|_{L^p(\Omega\times(T_0,T))}$ follows from the equation \eqref{main-eqn-2} and the Poincar\'e inequality. 
Therefore,  \eqref{global-time-cut-est} is proved by combining these three estimates.
\end{proof}

Finally, we provide the proof of Theorem \ref{main-theorem}.

\begin{proof}[Proof of Theorem \ref{main-theorem}]
By using the time dilation $(x,t)\mapsto (x,\lambda_0 t)$ as in \eqref{scale-2} and replacing $T$ by $T/\lambda_0$, we may assume without loss of generality that $\lambda_0=1$.
We first consider the case $p=2$ and begin by proving \eqref{main-thm-est}. Suppose that $u\in\W^{1,2}(\Omega_T,\beta)$ is a weak solution to \eqref{main-eqn}.  By Lemma \ref{embedd-lemma}, we see that $u \in L^2(\Omega_T, \beta)$. Hence,  by using Steklov's average as in \cite[p.~18]{DiB}, we can take $u$ as the test function of equation \eqref{main-eqn} to obtain
\[
\frac{1}{2}\frac{d}{dt}\int_{\Omega} \beta(x)u(x, t)^2  dx + \int_{\Omega}\wei{\bA \nabla u, \nabla u}dx = -\int_{\Omega} \wei{F, \nabla u}(x,t) dx+\int_{\partial \Omega}g(x,t)u(x,t)dS
\]
for a.e. $t\in (0, T)$. By \eqref{ellip-cond}, the Cauchy inequality, the trace theorem, and Young's inequality, it then follows that 
\begin{equation}\label{esti-p1}
\begin{aligned}
&\frac{d}{dt}\int_{\Omega} \beta(x)u(x, t)^2  dx +
\nu \|\nabla u(\cdot, t)\|_{L^2(\Omega)}^2\\
&\leq N(n,\nu,\Omega)\Big(\|u(\cdot, t)\|_{L^2(\Omega)}^2 +\|F(\cdot, t)\|_{L^2(\Omega)}^2+\|g(\cdot, t)\|_{L^2(\partial \Omega)}^2\Big),\quad \text{for a.e.}\quad t\in (0, T).
\end{aligned}
\end{equation}
To handle the term $\|u(\cdot, t)\|_{L^2(\Omega)}^2$, we recall that $\beta \in A_2$ with $[\beta]_{A_2}\leq M_0$ as $\frac{1}{\beta} \in A_{1+\frac{2}{n_0}}\subset A_{2}$ with $n_0=\max\{2,n\}$. Hence, by applying Lemma~\ref{embedd-lemma-2} with $\mu=\beta$ and using the fact that $\beta(\Omega)\leq |\Omega|^{1-\frac{2}{n_0}}\lambda_0=|\Omega|^{1-\frac{2}{n_0}}$ due to H\"older's inequality, we can infer that
\begin{equation*}
\|u(\cdot, t)\|_{L^2(\Omega)}^2\le \epsilon\|\nabla u(\cdot, t)\|_{L^2(\Omega)}^2+N(n, M_0, \Omega, \epsilon)\|u(\cdot, t)\|_{L^2(\Omega,\beta)}^2 \quad \text{for a.e.}~t\in(0,T).
\end{equation*}
Then, plugging this into \eqref{esti-p1} and choosing $\epsilon>0$ sufficiently small, we obtain
\begin{align*}
&\frac d{dt}\int_\Omega \beta(x) u(x,t)^2\,dx+\frac{\nu}{2}\|\nabla u(\cdot, t)\|_{L^2(\Omega)}^2\\
&\leq N(n,\nu,M_0,\Omega)\Big(\|u(\cdot, t)\|_{L^2(\Omega,\beta)}^2+\|F(\cdot, t)\|_{L^2(\Omega)}^2+\|g(\cdot,t)\|_{L^2(\partial\Omega)}^2\Big).
\end{align*}
By Gr\"onwall's inequality and integrating over $(0,T)$, it follows that
\begin{align*}
&\sup_{t\in(0,T)}\|u(\cdot,t)\|_{L^2(\Omega,\beta)}^2
+\int_0^T\|\nabla u(\cdot,t)\|_{L^2(\Omega)}^2\,dt\\
&\leq N\Big(\|u_0\|_{L^2(\Omega,\beta)}^2+\|F\|_{L^2(\Omega_T)}^2+\|g\|_{L^2(\partial\Omega\times(0,T))}^2\Big),
\end{align*}
where $N=N(n,\nu,M_0,\Omega, T)>0$.  Using Lemma~\ref{embedd-lemma-2} again, we obtain
\[
\|u\|_{L^2((0, T), W^{1,2}(\Omega))} \leq N \left(\|u_0\|_{L^2(\Omega, \beta)}+\|F\|_{L^2(\Omega_T)}+\|g\|_{L^2(\partial \Omega\times (0,T))}\right).
\]
Then, using the PDE in \eqref{main-eqn} and Lemma~\ref{embedd-lemma}, we infer that
\begin{equation} \label{main-p=2}
\|u\|_{\W^{1,2}(\Omega_T, \beta)} \leq N \left(\|u_0\|_{W^{1,2}(\Omega)}+\|F\|_{L^2(\Omega_T)}+\|g\|_{L^2(\partial \Omega\times (0,T))}\right).
\end{equation}
Hence, \eqref{main-thm-est} is proved. With this, the existence of a solution $u \in \W^{1,2}(\Omega_T, \beta)$ to \eqref{main-eqn} can be derived from the Galerkin method. The uniqueness of the solution also follows directly from \eqref{main-thm-est}.  Theorem \ref{main-theorem} is then proved when $p=2$.

\smallskip
Next, we consider the case $p >2$. In this case, we note that $F \in L^2(\Omega_T)^n$ and $g\in L^2(\partial \Omega\times (0,T))$ due to the assumptions on $F$ and $g$. Then, it follows from the case $p=2$ that there is a unique weak solution $u \in \W^{1,2}(\Omega_T, \beta)$ to \eqref{main-eqn}, and \eqref{main-p=2} holds. Thus, it remains to prove that $u \in \W^{1,p}(\Omega_T, \beta)$ and \eqref{main-thm-est} holds. To this end,  we trivially extend $u(\cdot, t)=u_0(\cdot)$ for all $t\in [-1, 0)$, $\bA(x,t)$ to be the identity matrix function on $\Omega\times (-1, 0)$,  $F=-\nabla u_0$ on $\Omega\times (-1, 0)$, and extend $g=0$ on $\partial \Omega\times (-1, 0)$. We note that $u \in \W^{1,2}(\Omega \times (-1, T), \beta)$, and it is a weak solution to 
\begin{equation} \label{main-eqn-extend}
\left\{
\begin{aligned}
\beta(x) u_t - \textup{div}(\bA \nabla u) & = \textup{div}(F) \quad  &&\text{in} \quad \Omega \times (-1, T), \\
(\bA \nabla u+F)\cdot \vec{\mathbf{n}}&=g \quad &&\text{on} \quad \partial \Omega\times (-1, T),\\ 
u(\cdot, -1) & = u_0(\cdot)  \quad &&\text{in} \quad \Omega.
\end{aligned} \right.
\end{equation}
Now, let $\bar\delta =\bar \delta(n, \nu, p, M_0) \in (0,1)$ be the number defined in Theorem~\ref{main-thm-2}. Note also that the conditions \eqref{ellip-cond}, \eqref{beta-cond}, and \eqref{small-ness} hold for the coefficients in \eqref{main-eqn-extend}. We then apply Theorem~\ref{main-thm-2}, the extension property above, and \eqref{main-p=2} to infer that $u \in \W^{1,p}(\Omega_T, \beta)$, and that
\begin{equation*}
   \begin{aligned}
\|u\|_{\W^{1,p}(\Omega_T,\beta)}
&\leq N\left( \|u\|_{L^2(\Omega\times (-1, T))} + \|F\|_{L^p(\Omega\times (-1, T))}+\|g\|_{L^{p}(\partial \Omega\times (-1, T))} \right)\\
&\leq N\left( \|u_0\|_{L^2(\Omega)}+\|u\|_{L^2(\Omega_T)}+\|\nabla u_0\|_{L^p(\Omega)} + \|F\|_{L^p(\Omega_T)}+\|g\|_{L^{p}(\partial \Omega\times (0, T))} \right)\\
&\leq N\left( \|u_0\|_{W^{1,p}(\Omega)}+\|F\|_{L^p(\Omega_T)}+\|g\|_{L^{p}(\partial \Omega\times (0, T))} \right),
\end{aligned}
\end{equation*}
where $N = N(n, \nu, p, M_0, R_0, \Omega, T) >0$. Thus, scaling back to the original time variable, we obtain \eqref{main-thm-est}. 

\smallskip
Lastly, we consider the case $p \in (1, 2)$. In this case, the existence, uniqueness, and the estimate \eqref{main-thm-est} can be obtained by the duality argument using the case $p \in (2, \infty)$. As the argument is standard, we skip the details. The proof of the theorem is completed.
\end{proof}

\section{Examples of Admissible Weights}\label{Examples}
We end the paper by providing some examples of admissible weights $\beta(x)$ satisfying \eqref{beta-cond} and \eqref{small-ness} in the following remark, whose detailed proof can be found in \cite[Appendix~A]{FP2}.

\begin{remark} \label{remark-example} 
The following examples satisfy the weight assumptions \eqref{beta-cond} and \eqref{small-ness} in
Theorem \ref{main-theorem}.
\begin{enumerate}[leftmargin=*]
\item[\textup{(a)}] Let $\beta_1(x) = |x|^\alpha$ with $x \in \mathbb{R}^n$. It follows from \cite{CMP, CFP} that $\frac{1}{\beta_1}(x)=|x|^{-\alpha} \in A_{1+\frac{2}{n_0}}$ when $\alpha \in (-2n/n_0, n)$. In addition, there exist $M_0 = M_0(n)\geq 1$ and $N = N(n) >0$ such that
\[
[\frac{1}{\beta_1}]_{A_{1+\frac{2}{n_0}}} \leq M_0 \quad \text{and} \quad [[\beta_1]]_{\textup{BMO}}(\Omega, R_0) \leq N|\alpha|, \quad \forall \ |\alpha| \leq \tfrac{1}{2}.
\]
Therefore, Theorem \ref{main-theorem} applies to $\beta_1(x)$ whenever $|\alpha|\in (0,\frac{1}{2})$ is sufficiently small. 

\smallskip
\item[\textup{(b)}]  Let $\beta_2: \mathbb{R}^n\setminus \{0\} \rightarrow \mathbb{R}$ be defined by
\begin{equation}\label{function varphi}
\beta_2(x) = \left\{
\begin{array}{ll}
-\ln|x| & \quad \textup{if} \quad |x| \leq e^{-1},\\
1& \quad \textup{otherwise}.
\end{array} \right.
\end{equation}
It is known that $\beta_2 \in A_1 \subset A_{1+\frac{2}{n_0}}$ with $[\beta_2]_{A_1} = M_0$ for some constant $M_0=M_0(n)\geq 1$; see \cite[Example~7.1.8]{Grafakos-2}, for instance. Moreover, there exists  $N = N(n)>0$ such that
\[
[[\beta_2]]_{\mathrm{BMO}}(B_{r_0},r_0/2)\leq\frac{N}{|\ln(4r_0)|}, \quad\text{for}\quad 0<r_0<1/(10e).
\]
Hence, $[[\beta_2]]_{\mathrm{BMO}}(B_{r_0},r_0/2)$ is sufficiently small if $r_0$ is sufficiently small. From this and a localization argument, we can apply Theorem \ref{main-theorem} to $\beta_2(x)$.

\smallskip
\item[\textup{(c)}] Given $m$ distinct points $\{x_k\}_{k=1}^m$ in $B_1$, let
\[
\omega_1(x) = \sum_{k=1}^{m} c_k |x-x_k|^{\alpha_k}, \quad \text{and} \quad \omega_2(x) = \prod_{k=1}^{m} |x-x_k|^{\alpha_k} \quad \text{for}\quad x \in \mathbb{R}^n,
\]
with $c_k>0$ and $\sum_{k=1}^m|\alpha_k|$ sufficiently small. Also, let
\begin{align*}
\omega_3(x)  =  \prod_{k=1}^{m}  \beta_2(x-x_k), \quad x \in \mathbb{R}^n,
\end{align*}
where $\beta_2$ is defined in \eqref{function varphi} of \textup{(b)}.  By using localization, dilation, translation, and a covering argument, we can apply Theorem \ref{main-theorem} to $\beta=\omega_i$ for $i\in \{1,2,3\}$.
\end{enumerate}
\end{remark}

\section*{Acknowledgement}
The research of T. Phan was partially supported by Simons Foundation, grant \# 769369.

\end{document}